\documentclass[letterpaper,11pt,reqno]{amsart} 
\usepackage[utf8]{inputenc}  
\usepackage[T1]{fontenc} 
\usepackage{indentfirst} 
\usepackage{amssymb,amsfonts} 
\usepackage{amsthm} 
\usepackage{amsmath} 
\usepackage[left=1in,right=1in,bottom=0.8in,top=0.8in]{geometry} 
\usepackage{multirow} 
\usepackage{enumerate} 
\usepackage{multicol} 
\usepackage{graphicx} 
\usepackage{setspace} 
\usepackage{color} 
\usepackage[all]{xy} 
\usepackage{fancyhdr} 
\usepackage{pdfpages} 
\usepackage{verbatim} 
\usepackage{float} 
\usepackage[labelfont=bf,format=plain]{caption} 
\usepackage{bbold} 
\usepackage{tikz}
\usepackage{mathtools}
\usepackage{relsize} 
\usepackage{url} 
\usepackage[hidelinks]{hyperref} 
\usepackage{subcaption}

\theoremstyle{plain} 
\newtheorem{teo}{Theorem}[]

\newtheorem{prop}{Proposition}[]
\newtheorem{lema}{Lemma}[] 
\newtheorem{cor}{Corollary}[]

\theoremstyle{definition} 

\newenvironment{dem}[1][Proof]{\noindent\textbf{#1:} }{\hfill \qed}

\newcommand{\R}{\mathbb{R}}

\newcommand{\Z}{\mathbb{Z}}

\newcommand{\PP}{\mathbb{P}}
\newcommand{\E}{\mathbb{E}}
\newcommand{\LL}{\mathcal{L}}
\newcommand{\EV}{\text{E}_{\text{vert}}}
\newcommand{\PPP}{\mathbb{P}_{p,q}^{\Lambda}}

\title[]{\bf Phase transition on randomly horizontally stretched square lattice}
\author[]{\small Isadora Guedes$^{*,\dagger}$, Paulo C Lima$^*$, Marcos Sá$^\diamond$, Remy Sanchis$^*$ }

\begin{document}
\pretolerance10000 
\setstretch{1.2} 

\begin{abstract} 
In this article, we study a bond percolation model on a horizontally stretched square lattice, constructed by stretching the distances between the columns of $\mathbb{Z}_+^2$ according to a collection of independent and identically distributed (i.i.d.) copies of a non-negative random variable $\xi$. We assume that $\xi$ satisfies the integrability condition
\[
\mathbb{E}\big[\xi\, e^{c(\log \xi)^{1/2}} \,\mathbb{1}_{\{\xi \geq 1\}}\big] < \infty,
\]
for some constant $c > 8\sqrt{\log 96}$. In this random environment, each vertical edge is independently declared open with probability $p$, while each horizontal edge is open with probability $p^{|e|}$, where $|e|$ denotes the Euclidean length of the edge. We develop a multiscale renormalization scheme adapted to this geometry and use it to prove that percolation occurs for all sufficiently large values of $p < 1$.
\end{abstract}
\maketitle





\maketitle

\section{Introduction}

Our goal in this work is to investigate the effect of introducing inhomogeneities on the lattice, specifically understanding how they affect the phase transition in percolation models. In some cases, inhomogeneities arise through the introduction of environments that specify the rules for assigning probabilities $p_e$ to each edge $e$ of the graph.

We will consider the square lattice, $\Z^2$, and one way of introducing inhomogeneities is by fixing certain columns, which will form what we will call the environment. In these columns, the edges will be open with probability $p$, while the remaining edges will be open with probability $q$, for some $p,q \in [0,1]$. Formally, given a subset $\Lambda \subset \Z$,  we define the set
\begin{equation*}
\EV = \big\{ \langle (x,y) , (x,y+1) \rangle : x \in \Lambda , y \in \Z \big\},
\end{equation*}
which contains the edges of those  columns which project into $\Lambda$. Given $p,q\in[0,1]$, an edge $e\in E(\Z^2)$ will be open with the probability
\begin{equation*}
p_e = \begin{cases} p, \text{ if } e \in \EV \\ q, \text{ if } e \notin \EV \end{cases}.
\end{equation*}
We will denote by $\PPP(\cdot)$ the quenched probability law in $\{0,1\}^{E(\Z^2)}$ which governs this percolation model and  we will refer to it  as percolation on the randomly horizontally stretched lattice in $\Z^2$.

In the special case where $\Lambda=\{0\}$, Zhang \cite{zhang} showed that $\PP_{p,q}^{\{0\}}((0,0) \leftrightarrow \infty) = 0$, for any $p \in [0,1)$ and $q \leq p_c(\Z^2) = \frac{1}{2}$. He used arguments similar to those of Harris in \cite{harris}, involving the construction of dual circuits around the origin, together with the Russo, Seymour and Welsh \cite{RSW1, RSW2} techniques. On the other extreme, when $\Lambda = \Z$, Kesten (see Section 11.9 in \cite{grim-book}) showed that $\PP_{p,q}^{\Z}((0,0) \leftrightarrow \infty) > 0$ if and only if  $p + q > 1$.
In the case where $\Lambda$ consists only of bounded gaps, namely, there exists $k \in \Z_+$ such that for every $l \in \Z, \; \Lambda \cap [l, l + k] \neq \emptyset$, a classic argument due to Aizenman and Grimmett \cite{aizenman-grimmett} ensured that for any $\epsilon >0$ there is $\delta = \delta(k, \epsilon) > 0$ such that $\PP_{p_c + \epsilon , p_c - \delta}^{\Lambda}((0,0) \leftrightarrow \infty) > 0$, where $p_c = p_c(\Z^2)$.

All the examples above involve deterministic environments.  We now turn to models in which  they  are random ones.

 For each $\rho \in [0, 1]$, let $\nu_{\rho}$ be the probability measure on subsets of $\Z$ in which the events $\{i \in \Lambda\}$ are independent, each occurring with probability $\rho$. Duminil-Copin et al \cite{brochette} showed that for any $\epsilon > 0$ and $\rho > 0$ there is a $\delta = \delta( \rho,\epsilon ) > 0$ such that $\PP_{p_c + \epsilon , p_c - \delta}^{\Lambda}((0,0) \leftrightarrow \infty) > 0$, for $\nu_{\rho}$-almost everywhere environment, where $p_c = p_c(\Z^2)$.
In \cite{BDS}, the authors considered the Contact Process on $\mathbb{Z}$ with an environment given by $\nu_{\rho}$. Their results, translated into Percolation, are that for any $\rho \in [0,1)$, there exists a sufficiently large $p<1$  such that $\PP_{0,p}(o \leftrightarrow \infty) > 0$ for $\nu_{\rho}$-almost every environment $\Lambda$. This means that if edges of columns are deleted according to Bernoulli trials with mean $\rho$, the phase transition is present. Hoffman \cite{hoffman} examined a similar case, where both rows and columns are independently deleted with the same probability $\rho$ and proved that the model still undergoes a non-trivial phase transition, see also \cite{marcos2}.

Hilário, Sá, Sanchis and Teixeira \cite{marcos} considered a random environment $\Lambda$ obtained by stretching horizontally the square lattice $\Z_+^2$ according to a positive random variable $\xi$, namely, they considered  independent and identically distributed copies of a positive random variable $\xi$ satisfying $\E(\xi^{\eta})< \infty$, for some $\eta>1$, to stretch the distance between the columns of $\Z_+^2$, obtaining a horizontally stretched square lattice.  By using a multiscale renormalization scheme,  they proved that the model exhibits a non-trivial phase transition. On the other hand, in  the same work, the authors also showed that, if $\mathbb{E}(\xi)=\infty$, the phase transition is trivial, leaving a gap on the moment conditions.

Aiming to fill  the gap mentioned above, we assume that $\E\big(\xi e^{c(\log \xi)^{1/2}} \mathbb{1}_{\{\xi \geq 1\}}\big) < \infty$, for a constant $c>8\sqrt{ \log 96}$, and using a finer multiscale renormalization scheme, our Theorem \ref{teo_trans_fase1},  shows that the model still undergoes a non-trivial phase transition.

In \cite{eulalia} the authors studied the so-called Renewal Contact Process on  $\mathbb{Z}^d$   and derived a condition—closely related to ours—ensuring the existence of a non-trivial subcritical phase. Despite of this similarity, both their renormalization scheme and their results differ from ours.

This article is organized as follows: in Section \ref{cap2} we describe the model and state Theorem \ref{teo_trans_fase1} which is our main result. In Section \ref{cap3} we prove the decoupling inequality given in Lemma \ref{lema_desig_desacop}. In Section \ref{cap4} we describe our multiscale scheme which consists of two parts: in Subsection \ref{secao_ambiente} we control the environments and in Subsection \ref{secao_cruzamentos} we control the vertical and the horizontal crossings. Finally, in Section \ref{secao_prova_teo}, we prove Theorem \ref{teo_trans_fase1}.

\section{The Model and Results}
\label{cap2}
Let $\Z_+$ be the set of all nonnegative integers and denote $\Z_+^* = \Z_+ \backslash \{0\}$. We will consider percolation in the lattice obtained from $\Z_+^2$ stretching randomly its horizontal edges. Formally speaking, let $\xi$ be a positive random variable and $\{\xi_i\}_{i\in\Z_+}$ a sequence of i.i.d. copies of $\xi$ and consider
\begin{equation*}
\Lambda = \left\{\sum_{1\leq i\leq k} \xi_i : k \in \Z_+ \right\},
\end{equation*}
which is called an environment. 
Notice that, $\Lambda$ can also be seen as the following increasing sequence 
\begin{equation}
\Lambda = \{x_k \in \R : x_0=0 \text{ and } x_k = x_{k-1} + \xi_k, \text{ for } k\in\Z_+^*\}.
\label{def_ambiente}
\end{equation}
This sequence is called a renewal process with interarrival distribution $\xi$. We will denote by $\mu_{\xi}(\cdot)$ the probability measure that governs this renewal process.

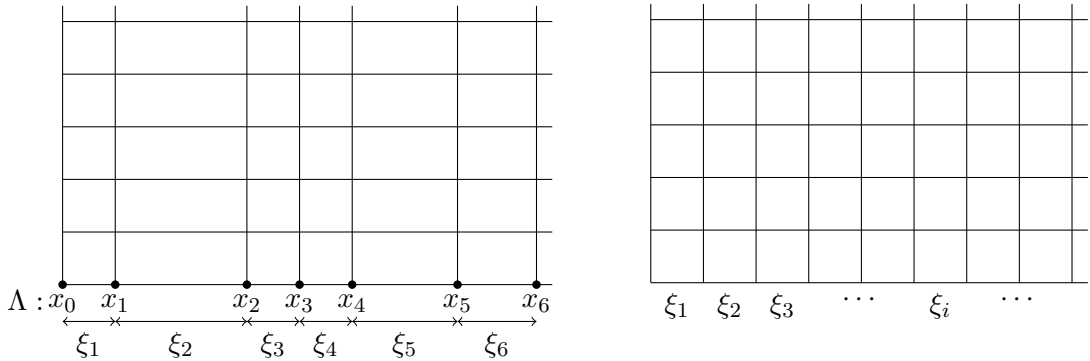
\begin{figure}[H]
     \centering
\subfloat[The lattice $\LL_{\Lambda}$, where the environment $\Lambda$ is given by $x_k$'s.]{
\begin{tikzpicture}[scale=.7]
\foreach \i in {0,1,3.5,4.5,5.5,7.5,9}{
\draw (\i,0)--(\i,5.3);
\filldraw[] (\i,0) circle (2pt);}
\foreach \j in {0,1,2,3,4,5}
\draw (0,\j)--(9.3,\j);
\node[below] at (-.7,0.09) {$\Lambda:$};
\node[below] at (0,0) {$x_0$};
\node[below] at (1,0) {$x_1$};
\node[below] at (3.5,0) {$x_2$};
\node[below] at (4.5,0) {$x_3$};
\node[below] at (5.5,0) {$x_4$};
\node[below] at (7.5,0) {$x_5$};
\node[below] at (9,0) {$x_6$};
\draw[<->] (0,-.7)--(1,-.7);
\node[below] at (0.5,-.7) {$\xi_1$};
\draw[<->] (3.5,-.7)--(1,-.7);
\node[below] at (2.25,-.7) {$\xi_2$};
\draw[<->] (3.5,-.7)--(4.5,-.7);
\node[below] at (4,-.7) {$\xi_3$};
\draw[<->] (5.5,-.7)--(4.5,-.7);
\node[below] at (5,-.7) {$\xi_4$};
\draw[<->] (5.5,-.7)--(7.5,-.7);
\node[below] at (6.5,-.7) {$\xi_5$};
\draw[<->] (9,-.7)--(7.5,-.7);
\node[below] at (8.25,-.7) {$\xi_6$};
\end{tikzpicture}
}
\hspace{1cm}
\subfloat[The lattice on $\Z_+^2$, where the environment $\Lambda$ is given by $\xi_k = x_k - x_{k-1}$.]{\begin{tikzpicture}[scale=.7]
\foreach \i in {0,1,...,8}
\draw (\i,0)--(\i,5.3);
\foreach \j in {0,1,2,3,4,5}
\draw (0,\j)--(8.3,\j);
\node[below] at (0.5,0) {$\xi_1$};
\node[below] at (1.5,0) {$\xi_2$};
\node[below] at (2.5,0) {$\xi_3$};
\node[below] at (4,0) {$\cdots$};
\node[below] at (5.5,0) {$\xi_i$};
\node[below] at (7,0) {$\cdots$};
\node[below] at (5,-1.225) {};
\end{tikzpicture}}
    \caption{Illustration of the lattice $\LL_{\Lambda}$ and its alternative formulation on $\Z_+^2$ with their respective environments.}
    \label{redes}
\end{figure}

Given a realization of an environment $\Lambda$ we can define the lattice $\LL_{\Lambda} = (V_{\Lambda}, E_{\Lambda})$ where the vertex and edge sets are given, respectively, by
\begin{equation*}
V_{\Lambda} = \Lambda \times \Z_+ = \{(x,y) \in \R^2 : x\in\Lambda, y\in\Z_+\}
\end{equation*}
and 
\begin{equation*}
E_{\Lambda} = \{\langle (x_i,n),(x_j,m) \rangle : |i-j| + |n-m| = 1, \text{ with } x_i, x_j \in \Lambda \text{ and } n,m \in \Z_+\},
\end{equation*}
see Figure \ref{redes}-(a).

Notice that  $\LL_{\Lambda}$ can be seen as the lattice obtained by $\Z_+^2$ by stretching or contracting the horizontal edges in such a way that $\xi_{i+1}$ gives the random separation between the $i$-th and $(i+1)$-th column in the stretched lattice.

We will also consider a bond percolation process in $\LL_{\Lambda}$ as follows. For each $p\in[0,1]$, denote by $\PP_p^{\Lambda}(\cdot)$ the probability measure on $\{0,1\}^{E_{\Lambda}}$ under which the random variables $\{\nu(e)\}_{e\in E_{\Lambda}}$ are independent Bernoulli random variables with mean
\begin{equation}
p_e = p^{|e|},
\label{peso_prob_elo_bernoulli} 
\end{equation}
where, for each $e=\langle v_1, v_2\rangle \in E_{\Lambda}$, $|e|=||v_1 - v_2||$ denotes the length of $e$, with $|| \cdot ||$ meaning the Euclidean norm in $\R^2$.
We say that an edge $e\in E_{\Lambda}$ is open if $\omega(e)=1$, and closed otherwise. We write $\{(0,0)\leftrightarrow \infty\}$ to represent the event that there is an infinite open path starting at $(0,0)$ which only uses open edges.

An equivalent formulation for the bond percolation model on $\LL_{\Lambda}$ defined by \eqref{peso_prob_elo_bernoulli} is the following: consider the square lattice $\Z_+^2$ and, conditional on $\xi_1, \; \xi_2, \dots$, declare each edge $e \in E(\Z_+^2)$ open independently with probability
\begin{equation}
p_e = \left\{\begin{array}{cc}
p, & \text{ if } e = \langle (i,j),(i,j+1) \rangle \text{ for some } i,j \\
p^{\xi_{i+1}}, & \text{ if } e = \langle (i,j),(i+1,j) \rangle \text{ for some } i,j
\end{array} \quad. \right. 
\label{peso_prob_elo_aberto}
\end{equation}

In the case that $\xi$ is a positive and integer-valued random variable, the percolation model defined on $\LL_{\Lambda}$ with parameters given by \eqref{peso_prob_elo_bernoulli} can still be mapped into another equivalent model on $\Z_+^2$ as follows. Consider the environment $\Lambda \subseteq \Z_+$ distributed as $\mu_{\xi}$ and define the set of edges
\begin{equation*}
E_{\text{vert}}(\Lambda^c) = \{\langle(x,y), (x,y+1) \rangle \in E(\Z_+^2) : x \notin \Lambda, \; y \in \Z_+\}.
\end{equation*}
Declare each edge $e\in\E(\Z_+^2)$ be open independently with probability 
\begin{equation}
p_e = \begin{cases}
0, \text{ if } e \in \E_{\text{vert}}(\Lambda^c) \\
p, \text{ if } e \notin E_{\text{vert}}(\Lambda^c)
\end{cases},
\label{peso_prob_elo_z2}
\end{equation}
and closed otherwise. This formulation is also a stretched square lattice obtained from $\Z_+^2$ by removing the edges lying in vertical columns that project to $\Lambda^c$ while preserving all other edges. Each one of these remaining edges is open independently with probability $p$. Notice that the resulting graph is similar to the stretched lattice $\LL_{\Lambda}$ defined above, however, the edges  now split into unit length segments, see Figure \ref{redes}-(b). We can recover the original formulation on $\LL_{\Lambda}$ by declaring an edge open if all the corresponding unitary edges in $\Z_+^2$ are open in the new formulation. 

Since all the formulations given above are equivalent, we will make a slight abuse of notation by also denoting the probability law of all versions by $\PP_p^{\Lambda}(\cdot)$.

Next we state our main result.

\begin{teo}
\label{teo_trans_fase1}
Let $\xi$ be a positive random variable that satisfies $\E \big(\xi e^{c(\log \xi)^{1/2}}\mathbb{1}_{\{\xi \geq 1\}} \big) < \infty$, for some constant $c>8\sqrt{\log 96}$. Then there is a critical point $p_c \in (0,1)$, depending on the law of $\xi$ only, such that for $p < p_c$, we have
\begin{equation}
\PP_p^{\Lambda}((0,0) \leftrightarrow \infty) = 0, \text{ for } \mu_{\xi}\text{-almost all } \Lambda,
\label{eq_ausencia_perc}
\end{equation}
and, for $p > p_c$, we have
\begin{equation}
\PP_p^{\Lambda}((0,0) \leftrightarrow \infty) > 0, \text{ for } \mu_{\xi}\text{-almost all } \Lambda.
\label{eq_existencia_perc}
\end{equation}
\end{teo}

The following sections will be devoted to the proof of Theorem \ref{teo_trans_fase1}.


\section{The Decoupling Inequality}
\label{cap3}
The goal of this section is prove the decoupling inequality given in Lemma \ref{lema_desig_desacop}, which will be essential in our multiscale scheme. We will begin with some definitions and notations about renewal processes. 

Let $\xi$ be a positive and $\chi$ a non-negative integer-valued random variable called, respectively, interarrival time and delay. Consider, as before, $\{\xi_i\}_{i\in\Z_+^*}$ i.i.d. copies of $\xi$ and also independent of $\chi$. We define recursively the renewal process $X = X(\xi, \chi) = \{X_i\}_{i\in\Z_+}$ as 
\begin{equation}
X_0=\chi, \text{  and  } X_i = X_{i-1} + \xi_i, \text{ for all } i \in \Z_+^*.  
\label{processo_renov_x}
\end{equation}

We say that the $i$-th renewal occurs at time $t$ if $X_{i-1}=t$. We will denote by $\mu_{\xi}^{\chi}(\cdot)$ the probability law that governs the renewal process $X$, regarded as a random element on a probability space supporting $\chi$ and the i.i.d. copies of $\xi$.

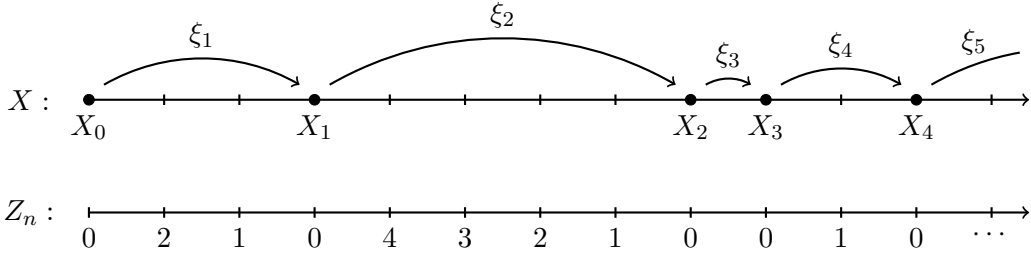
\begin{figure}[!htbp]
\centering
\begin{tikzpicture}[scale=1.0]
        
        \draw[thick,->] (0,1.5) -- (12.5,1.5);
        \foreach \x in {0,...,12}
        \draw[thick] (\x,1.43) -- (\x,1.57);
        \node at (-0.8,1.5) {$X:$};
        
        \foreach \x/\label in {0/X_0, 3/X_1, 8/X_2, 9/X_3, 11/X_4}
        {
            \filldraw[black] (\x,1.5) circle (2pt);
            \node[below=2pt] at (\x,1.5) {$\label$};
        }
        
        \draw[->, thick](0.2,1.7) arc[start angle=120, end angle=60, radius=2.6cm] node[midway, above] {$\xi_1$};
        \draw[->, thick] (3.2,1.7) arc[start angle=120, end angle=60, radius=4.6cm] node[midway, above] {$\xi_2$};
        \draw[->, thick] (8.2,1.7) arc[start angle=120, end angle=60, radius=0.6cm] node[midway, above] {$\xi_3$};
        \draw[->, thick] (9.2,1.7) arc[start angle=120, end angle=60, radius=1.6cm] node[midway, above] {$\xi_4$};
        \draw[-, thick] (11.2,1.7) arc[start angle=120, end angle=100, radius=3.6cm] node[midway, above] {$\xi_5$};
        
        \draw[thick,->] (0,0) -- (12.5,0);
        \foreach \x in {0,...,12}
        \draw[thick] (\x,-0.07) -- (\x,0.07);
        \node at (-0.8,0) {$Z_n:$};
        
        \foreach \x/\val in {0/0, 1/2, 2/1, 3/0, 4/4, 5/3, 6/2, 7/1, 8/0, 9/0, 10/1, 11/0, 12/\cdots}
        \node[below=2pt] at (\x,0) {$\val$};     
    \end{tikzpicture}
\caption{An illustration of the processes X and Z, when $\xi_1=3$, $\xi_2=5$, $\xi_3=1$ and $\xi_4=2$.}
\label{fig_processos_xyz}
\end{figure}

It is suitable to define other process related to $X$ (see Figure \ref{fig_processos_xyz}), namely
\begin{equation*}
Z = Z(\xi, \chi) = \{Z_n\}_{n\in\Z_+},
\end{equation*}
given by
\begin{equation}
Z_n = \min \{ X_i - n : i \in \Z_+ \text{ and } X_i - n \geq 0\}. 
\label{processo_renov_z}
\end{equation}

Notice that knowing each one of the processes $X$ or $Z $ we are able to determine the other. Thus, with an abuse of notation, we will also refer to the process $Z$ as renewal processes with interarrival time $\xi$ and delay time $\chi$.  Additionally, we will use $\mu_{\xi}^{\chi}(\cdot)$ to denote the probability law of the renewal process $Z$.

It is important to notice that $Z$ is a Markov chain with transition probability given by
\begin{equation}
\PP(Z_n = i | Z_{n-1} = j) = \begin{cases}
\PP(\xi = i+1), \text{ if } j=0 \\
1, \text{ if } j = i+1  \\
0, \text{ otherwise}
\end{cases},
\label{prob_transicao}
\end{equation}
for all $n\in\Z_+^*$ and $i,j \in \Z_+$.

When $\E(\xi)< \infty$, we can define a random variable $\rho = \rho(\xi)$ with distribution 
\begin{equation}
\lambda_k = \PP(\rho = k) := \dfrac{1}{\E(\xi)} \sum_{i=k+1}^{\infty} \PP(\xi = i), \text{ for all } k\in\Z_+,
\label{medida_estacionaria}
\end{equation}
independent of everything else. Consider the renewal process $Z(\xi, \rho)$, defined in \eqref{processo_renov_z}, using $\rho$ as its delay time. We can show by induction that 
\begin{equation*}
Z_n \stackrel{d}{=} Z_0 \stackrel{d}{=} \rho, \text{ for all } k \in \Z_+.
\end{equation*}
For this and since $Z$ is a Markov chain, we have that
\begin{equation}
\theta_m Z \stackrel{d}{=} Z, \text{ for any } m \in \Z_+^*,
\label{inv_shift}
\end{equation}
where $\theta_m$ is the left shift operator, namely, $\theta_m: \Z^{\infty} \rightarrow \Z^{\infty}$, given by
\begin{equation*}
\theta_m(x_0, x_1, \dots) = (x_m, x_{m+1}, \dots).
\end{equation*}
For this reason, for a fixed $\xi$, the random variable $\rho$, given in \eqref{medida_estacionaria}, is called stationary delay.

For each $c >0$, let
\begin{equation}
F_c(x)= e^{c(\log x)^{1/2}},  
\label{def_Fc}
\end{equation}
where $x \geq 1$. Using \eqref{medida_estacionaria}, we can show the following result.

\begin{lema}
\label{lema_momento_dist_estacionaria}
If $\xi$ is a positive random variable satisfying $\E\big(\xi F_c(\xi) \mathbb{1}_{\{\xi \geq 1\}}\big) < \infty$ and $c>0$, then the stationary delay $\rho = \rho(\xi)$, given by \eqref{medida_estacionaria}, satisfies $\E \big(F_c(\rho) \mathbb{1}_{\{\rho \geq 1\}}\big)< \infty$.
\end{lema}

We say that a random variable $\xi$ is aperiodic if
\begin{equation*}
\gcd \{ k \in \Z_+^* : \PP(\xi=k) >0\} =1.
\end{equation*}

In what follows we will assume that $\xi$ is aperiodic and  $\E(\xi)<\infty$.  Let $Z=Z(\xi, \chi)$ and $Z'=Z'(\xi,\chi')$ be two independent renewal processes with interarrival time $\xi$ and delays $\chi$ and $\chi'$, respectively. Define
\begin{equation}
T = \min\{k\in\Z_+^* : Z_k = Z'_k=0\},
\label{def_T}
\end{equation}
as the coupling time of $X$ and $X'$. We will denote by $\nu^{\chi, \chi'}_{\xi}(\cdot)$ the product measure $\nu_{\xi}^{\chi} \otimes \nu_{\xi}^{\chi'}$.

Now we will establish the decoupling inequality for stationary renewal processes given in Lemma \ref{lema_desig_desacop}. For that, using the Markov inequality, we will need to get an upper bound of the form $\mu_{\xi}^{\chi, \chi'}(T > n) \leq \dfrac{c_0}{F_c(n)}$, where $c_0$ is a constant. Define
\begin{equation}
\widetilde{c}=\max\Big\{e^{\left(\frac{\log 2}{c}\right)^2}, e^{1/2}\Big\}.
\label{cota_c_teo}
\end{equation}

\begin{lema}[Proposition 3, \cite{lindvall}]
\label{teo_processo_renov}
Suppose that $\xi$ is an aperiodic and integer-valued random variable taking values greater that $\widetilde{c}$ and satisfies $\E \big(\xi F_c(\xi) \big) < \infty$, where $F_c$ is given in \eqref{def_Fc}. Also suppose that $\chi$, $\chi'$ are non-negative integer-valued random variables with $\E \big(F_c(\chi) \mathbb{1}_{\{\chi \geq 1\}} \big)$ and $\E \big(F_c(\chi') \mathbb{1}_{\{\chi' \geq 1\}} \big)$ finite. Then, $\E_{\xi}^{\chi, \chi'} \big(F_c(T)\mathbb{1}_{\{T \geq 1\}} \big) < \infty$, where $T$ is given by \eqref{def_T}.
\end{lema}

\begin{dem}
From Proposition 3 of \cite{lindvall}, in order to prove the lemma above it is enough to show that the function $F_c(x)$ satisfies the following conditions (which appear on page 63 of \cite{lindvall}):
\begin{enumerate}[(i)]
    \item $F_c(x)$ is non-decreasing and $F_c(x) \geq 2$, for all $x \geq \widetilde{c}$; 
    \item $\dfrac{\log F_c(x)}{x}$ is non-increasing, for all $x \geq \widetilde{c}$ and $\displaystyle\lim_{x \rightarrow \infty} \frac{\log F_c(x)}{x} = 0$.
\end{enumerate}

Notice that conditions (i) and (ii) follow from the definition of $\widetilde{c}$, which concludes the proof.
\end{dem}

The next lemma gives the desired decoupling inequality for stationary renewals. 

\begin{lema}[{\bf Decoupling Inequality}]
\label{lema_desig_desacop}
Let $\xi$ be an aperiodic and integer-valued random variable taking values greater than $\widetilde{c}$ and satisfies $\E \big(\xi F_c(\xi) \big) < \infty$,  where $F_c$ is given in \eqref{def_Fc}. Consider the associated renewal process $Z$ defined in \eqref{processo_renov_z}. Then, there exists a constant $c_1 = c_1(\xi, c) > 0$ such that for all $n,m\in\Z_+$ and for every pair of events $A$ and $B$, with
\begin{equation*}
A \in \sigma(Z_j : 0 \leq j \leq m) \text{  and  } B\in\sigma(Z_j : j \geq m+2n),
\end{equation*}
we have
\begin{equation}
\mu_{\xi}^{\rho}(A \cap B) \leq \mu_{\xi}^{\rho}(A) \mu_{\xi}^{\rho}(B) + \dfrac{c_1}{F_c(n)}.
\label{desig_descop}
\end{equation}
\end{lema}

\begin{dem}
If $\mu_{\xi}^{\rho}(A)=0$, we are done. So let us suppose that $\mu_{\xi}^{\rho}(A)>0$. Using the definition of the renewal process $Z$ given in \eqref{processo_renov_z}, we have

\begin{align}
\nonumber \mu_{\xi}^{\rho}(A \cap B) & = \mu_{\xi}^{\rho}(A \cap B \cap \{Z_m > n\}) + \mu_{\xi}^{\rho}(A \cap B \cap \{Z_m \leq n\}) \\
\nonumber & \leq \mu_{\xi}^{\rho}(Z_m > n) + \mu_{\xi}^{\rho}(A) \mu_{\xi}^{\rho}(B \cap \{Z_m \leq n\}|A) \\
\nonumber & = \mu_{\xi}^{\rho}(Z_m > n) + \mu_{\xi}^{\rho}(A) \cdot \sum_{\substack{0 \leq i \leq n \\ \mu_{\xi}^{\rho}(Z_m=i|A)>0}} \mu_{\xi}^{\rho}(B |A \cap \{Z_m =i\}) \mu_{\xi}^{\rho}(Z_m =i |A) \\ 
\nonumber & = \mu_{\xi}^{\rho}(Z_m > n) + \mu_{\xi}^{\rho}(A) \cdot \sum_{\substack{0 \leq i \leq n \\ \mu_{\xi}^{\rho}(Z_m=i|A)>0}} \mu_{\xi}^{\rho}(B |\{Z_m =i\}) \mu_{\xi}^{\rho}(Z_m =i |A) \\
\nonumber & = \mu_{\xi}^{\rho}(Z_m > n) + \mu_{\xi}^{\rho}(A) \max_{0 \leq i \leq n} \mu_{\xi}^{\delta_{m+i}}(B) \cdot \sum_{\substack{0 \leq i \leq n \\ \mu_{\xi}^{\rho}(Z_m=i|A)>0}} \mu_{\xi}^{\rho}(Z_m =i |A) \\
 & \leq \mu_{\xi}^{\rho}(Z_m > n) + \mu_{\xi}^{\rho}(A) \max_{0 \leq i \leq n} \mu_{\xi}^{\delta_{m+i}}(B).
\label{des_desacop}
\end{align}

Now, we need compare the measures $\mu_{\xi}^{\delta_{m+i}}$ and $\mu_{\xi}^{\rho}$, when $0\leq i \leq n$. For this, notice that $\mu_{\xi}^{\delta{m+i}}(B) = \mu_{\xi}^{\delta_0}(\theta_{m+i}(B))$ and that we have given a interval of size $2n-i$ for the renewal processes with delays $\delta_0$ and $\rho$ to couple. So, by the stationarity of $\rho$, we have
\begin{align*}
|\mu_{\xi}^{\delta_{m+i}}(B) - \mu_{\xi}^{\rho}(B)| & = |\mu_{\xi}^{\delta_0}(\theta_{m+i}(B)) - \mu_{\xi}^{\rho}(\theta_{m+i}(B))| \\
 & = \big| \mu_{\xi}^{\delta_0} \big(\theta_{2n-i}Y \in \theta_{m+i}(B) \big) - \mu_{\xi}^{\rho} \big(\theta_{2n-i}Y \in \theta_{m+i}(B) \big) \big| \\
 & \leq \mu_{\xi}^{\delta_0,\rho} (T > 2n-i) \\
 & \leq \mu_{\xi}^{\delta_0,\rho} (T > n),
\end{align*}
where the last inequality follows since $0 \leq i \leq n$. Thus, for all $0\leq i \leq n$, we get
\begin{equation*}
\mu_{\xi}^{\delta_{m+i}}(B) \leq \mu_{\xi}^{\rho}(B) + \mu_{\xi}^{\delta_0,\rho} (T > n).
\end{equation*}
Using this inequality in \eqref{des_desacop} and using that $Z_m \stackrel{d}{=} Z_0  \stackrel{d}{=} \rho$, we have

\begin{align*}
\mu_{\xi}^{\rho}(A \cap B) & \leq \mu_{\xi}^{\rho}(\rho > n) + \mu_{\xi}^{\rho}(A) \mu_{\xi}^{\rho}(B) + \mu_{\xi}^{\rho}(A) \mu_{\xi}^{\delta_0,\rho} (T > n) \\
 & \leq \mu_{\xi}^{\rho}(A) \mu_{\xi}^{\rho}(B) + \mu_{\xi}^{\rho}(\rho > n) + \mu_{\xi}^{\delta_0,\rho} (T > n) \\
 & \leq \mu_{\xi}^{\rho}(A) \mu_{\xi}^{\rho}(B) + \dfrac{\E_{\xi}^{\rho}\big(F_c(\rho)\big)}{F_c(n)} + \frac{\E_{\xi}^{\delta_0, \rho}\big(F_c(T)\big)}{F_c(n)},
\end{align*}
where the last inequality follows from the Markov inequality for $F_c(\rho)$ and $F_c(T)$. We conclude the proof by taking $c_1 = \E_{\xi}^{\rho}\big(F_c(\rho)\big) + \E_{\xi}^{\delta_0, \rho}\big(F_c(T)\big)$, which is finite by Lemmas \ref{lema_momento_dist_estacionaria} and \ref{teo_processo_renov}.
\end{dem}

\section{The Multiscale Scheme}
\label{cap4}

In this section, we will present our multiscale scheme, which consists of two main parts: environments and crossing events. In Subsection \ref{secao_ambiente}, we will define an increasing sequence of numbers $L_k$, $k\geq 0$, called scales, and use them to partition the set $\mathbb{Z}_+$ into intervals of length $L_k$, named $k$-intervals. We will find some integer $k_0$, which will satisfy some conditions that will appear throughout the text, and we will label the $k$-intervals as good or bad, recursively, for all $k \geq k_0$. Next, we will show that, with strictly positive probability, there is an environment that has the following property: the first two $k$-intervals are good, for all $k \geq k_0$. In Subsection \ref{secao_cruzamentos}, we will construct horizontal and vertical crossings of rectangles whose sides depend on $L_k$ and prove that, within $k$-good intervals, such crossings have a very high probability to occur. This allows us to build an infinite open cluster in the proof of Theorem \ref{teo_trans_fase1}.

\subsection{Environments}\mbox{}\\
\label{secao_ambiente}

Given $c>0$, for any $e^{-c^2/32}<\alpha<1$, let
\begin{equation}
A = A(\alpha,c) = \alpha e^{c^2/32}.
\label{def_A}
\end{equation}
Notice that $1<A<e^{c^2/32}$, in particular,  \begin{equation}\label{eqrest}\dfrac{c}{4}>\sqrt{2\log A}.\end{equation}

Consider the sequence of numbers $(L_k)_{k \in \Z_+}$, called scales, defined recursively by
\begin{equation}
L_0 = \lfloor A \rfloor \quad \text{and} \quad L_{k} = \big\lfloor A^{k+1} \big\rfloor L_{k-1}, \text{ for all } k \geq 1.
\label{def_escalas}
\end{equation}

Given $x>1$, then $\dfrac{x}{2} \leq \lfloor x\rfloor \leq x$. In fact, the upper bound follows by definition of the function $\lfloor \cdot \rfloor$. To get the lower bound, notice that if $1<x<2$, then $\lfloor x \rfloor=1>\dfrac{x}{2}$. On the other hand, if $x \geq 2$,  then   $\lfloor x \rfloor\geq x-1=x(1-1/x)\geq \dfrac{x}{2}$. Therefore,  for all $k\geq 0$, we have
\begin{equation}
\dfrac{A^{k+1}}{2} \leq \lfloor A^{k+1} \rfloor \leq A^{k+1}.
\label{cotas_chao_Ak}
\end{equation}

From  \eqref{def_escalas} and \eqref{cotas_chao_Ak}, it follows by induction,  that  for $k\geq 1$, we have
\begin{equation}
\dfrac{A^{\frac{(k+1)(k+2)}{2}}}{2^{k+1}} \leq L_{k} \leq A^{\frac{(k+1)(k+2)}{2}},
\label{cotas_Lk}
\end{equation}
which implies 
\begin{equation}
L_k \leq A^{\frac{(k+1)(k+2)}{2}} \leq 2^{k+1} L_k,
\label{cotas_Ax}
\end{equation}
and so,
\begin{equation}
L_k^{\frac{2}{k+1}} \leq A^{k+2} \leq 4\, L_k^{\frac{2}{k+1}}.
\label{qqc}
\end{equation}
Using (\ref{qqc}),  \eqref{def_escalas} and \eqref{cotas_chao_Ak}, we have
\begin{equation}
\dfrac{1}{2}\,  L_k^{\frac{k+3}{k+1}} \leq L_{k+1} \leq 4\,  L_k^{\frac{k+3}{k+1}}
\label{cotas_Lk+1_e_Lk}
\end{equation}
and from \eqref{cotas_Lk+1_e_Lk} and \eqref{cotas_Lk}, respectively,  it follows that
\begin{align}
\log L_{k+1} & \leq \left(\dfrac{\log 4}{\log L_k} + \dfrac{k+3}{k+1}\right) \log L_k \nonumber \\
& \leq  \left(\dfrac{2\log4}{(k+1)\big[(k+2) \log A-2\log2\big]} + \dfrac{k+3}{k+1}\right) \log L_k \nonumber \\
& =: r(k,A) \, \log L_k.
\label{eqsep24}
\end{align}
Also, from \eqref{cotas_Lk}, we obtain
\begin{equation}
\dfrac{4}{k+1} (\log L_k)^{1/2} \leq 2\sqrt{2} \left(\dfrac{k+2}{k+1}\right)^{1/2}(\log A)^{1/2}.
\label{eq1}
\end{equation}

Let $c_1$ be the constant given in  Lemma \ref{lema_desig_desacop} and take $k_0=k_0(c, \alpha, c_1) \in \Z_+$ sufficiently large such that for all $k \geq k_0$ the following three conditions (which will be necessary in the  prove of Lemma \ref{lema_cota_pk_lindvall}) hold:
\begin{equation}
c -2 \left(\dfrac{k+2}{k+1}\right)^{1/2}(2\log A)^{1/2} - \dfrac{c}{2} \big(r(k,A)\big)^{1/2} >  \dfrac{c}{4} - \sqrt{2\log A},
\label{cond_k0_1}
\end{equation}

\begin{equation}
\exp\left[-\left(\dfrac{c}{4} - \sqrt{2 \log A} \right) \left(\dfrac{(k+1)(k+2)}{2}\right)^{1/2}(\log A)^{1/2}\right] < \dfrac{1}{16(c_1 + 1)},
\label{cond_k0_2}
\end{equation}

\begin{equation}
\label{cond_k0_3}
\E\big(e^{c(\log \rho)^{1/2}}\mathbb{1}_{\{\rho\geq 1\}}\big)\leq \exp\left(\dfrac{c\,\,(\log L_{k})^{1/2}}{2}\right) .
\end{equation}

For each $k \geq k_0$, we divide $\Z_+$ into disjoint intervals of length $L_k$. The $i$-th interval of scale $k$, denoted by $I_i^k$, is defined by
\begin{equation}
I_i^k = [iL_k, (i+1)L_k), \; \text{for } i\in\Z_+.
\label{intervalos}
\end{equation}
Notice that each interval of scale $k$ can be partitioned into $ \lfloor A^{k+1} \rfloor$ sub-intervals of scale $k-1$, namely, 
$$I_i^k = \bigcup_{l\in\mathcal{I}_{k,i}} I_l^{k-1},$$
where $\mathcal{I}_{k,i} = \left\{i \lfloor A^{k+1} \rfloor , i \lfloor A^{k+1} \rfloor  +1, \dots, (i+1) \lfloor A^{k+1} \rfloor  -1 \right\}$ represents the set of indices for the sub-intervals of scale $k-1$ which are within of $I_i^k$.

Given an environment $\Lambda \subseteq \Z_+$ as in \eqref{def_ambiente} and $k \geq k_0$, we will label the intervals $I_i^k$ either as {\it good} or {\it bad}, recursively. For $k=k_0$, we declare $I_i^{k_0} \;$ good if $\Lambda \cap I_i^{k_0} \neq \emptyset$, that is, if there exists at least one column of $\LL_{\Lambda}$ present; otherwise, we declare it bad. For $k > k_0$, assuming that all intervals at scale $k-1$ have been defined, we declare  $I_i^k \;$ bad if it has at least two non-consecutive bad intervals of the scale $k-1$; otherwise, we declare it as good. Notice that a good interval at scale $k$ can have a maximum of two bad intervals of the scale $k-1$ and, in this case, these intervals must be adjacent because otherwise we would have two non-consecutive bad intervals.

For each $i\in\Z_+$ and $k\geq k_0$, let $A_i^k$ be the event
$$A_i^k=\{I_i^k \text{ is bad}\}.$$
Sometimes we will write $\{I_i^k \text{ is good}\}$ for the complement of $A_i^k$. Define
\begin{equation}
p_k := \mu_{\xi}^{\rho}(A_0^k) = \mu_{\xi}^{\rho}(A_i^k),
\label{definicao_pk}
\end{equation}
where the last equality follows from the stationarity of $\rho$.

To get an upper bound for $p_{k+1}$ in terms of $p_k$, we notice that
\begin{align}
\nonumber p_{k+1} & = \mu_{\xi}^{\rho}(\text{there is at least two non-consecutive bad intervals on the scale } k) \\ 
\nonumber & \leq \binom{\lfloor A^{(k+2)} \rfloor}{2} \mu_{\xi}^{\rho}\big(A_0^k \cap A_2^k\big) \\ 
\nonumber & \leq \binom{\lfloor A^{(k+2)} \rfloor}{2} \left(p_k^2 + \dfrac{c_1}{e^{c(\log L_k)^{1/2}}} \right) \\
& \leq \frac{A^{2(k+2)}}{2} \left(p_k^2 + \dfrac{c_1}{e^{c(\log L_k)^{1/2}}} \right),
\label{relacao_pk}
\end{align}
where the second inequality follows from  the decoupling inequality given in Lemma \ref{lema_desig_desacop}.

The next result shows that the probability of an interval being bad decreases exponentially fast in $k$, for all scales $k\geq k_0$. 

\begin{lema}
\label{lema_cota_pk_lindvall}
Let $c>0$ and $k_0$ satisfying \eqref{cond_k0_1}-\eqref{cond_k0_3}, then for all $k \geq k_0$, we have
\begin{equation}
p_k \leq \exp{\left(-\dfrac{c(\log L_k)^{1/2}}{2}\right)}.
\label{pk}
\end{equation}
\end{lema}

\begin{dem}
The proof of this lemma will be done by induction in $k$. For $k=k_0$, using \eqref{definicao_pk}, the stationarity of $\rho$, the Markov inequality and \eqref{cond_k0_3}, respectively, we get
\begin{equation*}
p_{k_0} = \mu_{\xi}^{\rho}(A_0^{k_0}) = \mu_{\xi}^{\rho}(Z_0 > L_{k_0}) = \PP(\rho > L_{k_0}) \leq \dfrac{\E\big(e^{c(\log \rho)^{1/2}} \mathbb{1}_{\{\rho\geq 1\}}\big)}{e^{c(\log L_{k_0})^{1/2}}} \leq \exp\left(-\dfrac{c(\log L_{k_0})^{1/2}}{2}\right).
\end{equation*}
Now suppose that, for some $k \geq k_0$, we have
\begin{equation}
p_k \leq \exp\left(-\dfrac{c(\log L_k)^{1/2}}{2}\right),
\label{hipotese_pk}
\end{equation}
then  
\begin{align}
p_{k+1} & \leq \dfrac{A^{2(k+2)}}{2} \left(p_k^2 + \dfrac{c_1}{e^{c(\log L_k)^{1/2}}} \right) \nonumber \\
 & \leq \dfrac{A^{2(k+2)}}{2} \left(e^{-c(\log L_k)^{1/2}} + c_1 e^{-c(\log L_k)^{1/2}} \right) \nonumber\\
 & = \dfrac{c_1 +1}{2} \; A^{2(k+2)} e^{-c(\log L_k)^{1/2}} \nonumber\\
 & \leq \dfrac{c_1 +1}{2} \; 4^2 L_k^{\frac{4}{k+1}} e^{-c(\log L_k)^{1/2}}, \label{cotapkx}
\end{align}
where the inequalities above follow from \eqref{relacao_pk}, \eqref{hipotese_pk} and \eqref{qqc} respectively. 
Hence, using  \eqref{cotapkx},   \eqref{eqsep24}, \eqref{eq1}, \eqref{cond_k0_1}, \eqref{cotas_Lk}  and \eqref{cond_k0_2}, respectively,  we get
\begin{eqnarray*}
\mbox{}&\mbox{} & \dfrac{p_{k+1}}{\exp\left(-\frac{c(\log L_{k+1})^{1/2}}{2}\right)}\nonumber\\
& \leq &  \dfrac{c_1+1}{2} \, 16 \, L_k^{\frac{4}{k+1}} e^{-c(\log L_k)^{1/2}} \exp\left(\frac{c(\log L_{k+1})^{1/2}}{2}\right)\nonumber \\
& = & 8 (c_1+1) \; \exp\left[\dfrac{4}{k+1}\log L_k - c(\log L_k)^{1/2} + \dfrac{c(\log L_{k+1})^{1/2}}{2}\right] \\
& \leq  & 8 (c_1+1) \; \exp\left[-(\log L_k)^{1/2} \left(c - \dfrac{4}{k+1} \left(\dfrac{(k+1)(k+2)}{2}\right)^{1/2}(\log A)^{1/2}   - \dfrac{c}{2} \big(r(k,A)\big)^{1/2}\right)\right] \\
   & \leq & 8 (c_1+1)\exp\left[- \left(\dfrac{c}{4} - \sqrt{2 \log A} \right) (\log L_k)^{1/2} \right] \\
& \leq  & 8 (c_1+1) \exp\left[- \left(\dfrac{c}{4} - \sqrt{2 \log A} \right) \left(\dfrac{(k+1)(k+2)}{2}\right)^{1/2}(\log A)^{1/2} \right] \label{cota_p_k+1_limite_exp}  \\
 & \leq  & 8 (c_1+1) \; \dfrac{1}{16(c_1+1)}\\
 &<&1,
 \end{eqnarray*}
 which concludes the proof of the lemma.
\end{dem}

Using the upper bound of $p_k$ that we have just proved, it follows that with strictly positive probability, the environment $\Lambda$ has the property that the intervals $I_0^k$ and $I_1^k$ are good for all scale $k$, with $k\geq k_0$. Moreover, we can show that all their subintervals are also good.

\begin{cor} Consider $I_i^k$ given by \eqref{intervalos} and let $k_0 \in \Z_+$ satisfying the conditions \eqref{cond_k0_1}-\eqref{cond_k0_3}. Then,
\begin{equation}
\mu_{\xi}^{\rho} \left(\bigcap_{k\geq k_0} \{I_0^k \text{ and } I_1^k \text{ are good}\}\right)>0.
\label{prob_blocos_bons}
\end{equation}
Moreover, 
\begin{equation}
\mu_{\xi}^{\rho} \left(\bigcap_{k\geq k_0} \bigcap_{i=0}^{8L_k^{\frac{2}{k+1}}-1} \{I_i^k \text{ is good}\}\right)>0.
\label{prob_blocos_bons_new}
\end{equation}
\end{cor}

\begin{dem}Observe that
   \begin{eqnarray*} \mu_{\xi}^{\rho} \left(\bigcap_{i=0}^{8L_k^{\frac{2}{k+1}}-1} \{I_i^k \text{ is good}\}\right)&\leq& 8L_k^{\frac{2}{k+1}}\exp{\left(-\dfrac{c(\log L_k)^{1/2}}{2}\right)}\\
   &\leq& \exp{\left(-\dfrac{c(\log L_k)^{1/2}}{2}+\log 8 + \dfrac{2}{k+1}\log L_k\right)}
  \end{eqnarray*} 
   is summable in $k$, since that $c$ is large. Equation \eqref{prob_blocos_bons_new} follows from Borel–Cantelli.
\end{dem}

\subsection{Crossing Events}\mbox{}\\
\label{secao_cruzamentos}

In this section we will study the probability of crossing events within special rectangles in $\Z^2_+$. The bases of these rectangles will be intervals on some scale $k$ and their heights will be much greater than their bases. This elongated form will be important so that we have many chances to cross these rectangles horizontally. Then we will use these rectangle crossings to build an infinite cluster.

Before stating our results, we need to introduce some notations and define some crossing events. Let $a,b,c,d\in\Z_+$ with $a<b$ and $c<d$, let $[a,b]=\{i \in \Z_+ : a \leq i \leq b\}$ and $[c,d]=\{i \in \Z_+ : c \leq i \leq d\}$. We define the rectangle $R$, denoted by
\begin{equation}
R=R\big([a,b)\times[c,d)\big)
\label{retangulo}
\end{equation}
as the subgraph of $\Z^2_+$ whose vertex and edge sets are given, respectively, by
\begin{equation*}
V(R)=[a,b]\times[c,d]
\end{equation*}
and
\begin{equation*}
E(R)=\big\{\langle(x,y),(x+l,y+1-l)\rangle \; : \; (x,y)\in[a,b-1]\times[c,d-1], \; l\in\{0,1\}\big\}.
\end{equation*}
In other words, $R$ is the rectangle $[a,b]\times[c,d]$ with the edges in the top and right sides removed.

We define the horizontal and vertical crossing events in $R$, denoted, respectively, by $\mathcal{C}_h(R)$ and $\mathcal{C}_v(R)$, as
\begin{equation}
\mathcal{C}_h(R)=\big\{\{a\}\times[c,d]\leftrightarrow \{b\} \times[c,d] \text{ in } R\big\}
\label{cruzamento_horiz}
\end{equation}
and
\begin{equation}
\mathcal{C}_v(R)=\big\{[a,b]\times \{c\} \leftrightarrow[a,b]\times \{d\} \text{ in } R\big\}.
\label{cruzamento_vert}
\end{equation}

Let us proceed by defining the specific rectangles and crossings that are of interest to us. Consider parameters $\beta, \mu \in (0,1)$ with $\mu < \beta$ and define a sequence of heights $H_0,H_1, \dots,$ recursively as follow: 
\begin{equation}
H_0 = 100 \text{ and } H_k = 2 \; \left \lceil e^{L_k^{\left(1-\frac{\beta}{k+1}\right)}} \right \rceil H_{k-1}, \; \text{ for all } k \geq 1.
\label{def_alturas}
\end{equation}
The choice of the initial height $H_0$ is really arbitrary and we could have chosen any other arbitrary positive integer.

We will consider 2 types of rectangles: one whose base is formed by 2 consecutive intervals on some scale $k$ and whose height is $H_k$ and another whose base consists of one interval on some scale $k$ and whose height is $2H_k$. For each $i,j,k\in\Z_+$, we will denote the horizontal and vertical crossings, at scale $k$, respectively, as follows
\begin{equation}
H_{i,j}^k = \mathcal{C}_h\left((I_i^k \cup I_{i+1}^k) \times [jH_k,(j+1)H_k)\right)
\label{evento_cruz_horiz}
\end{equation}
and
\begin{equation}
V_{i,j}^k =\mathcal{C}_v\left(I_i^k \times [jH_k,(j+2)H_k)\right).
\label{evento_cruz_vert}
\end{equation}
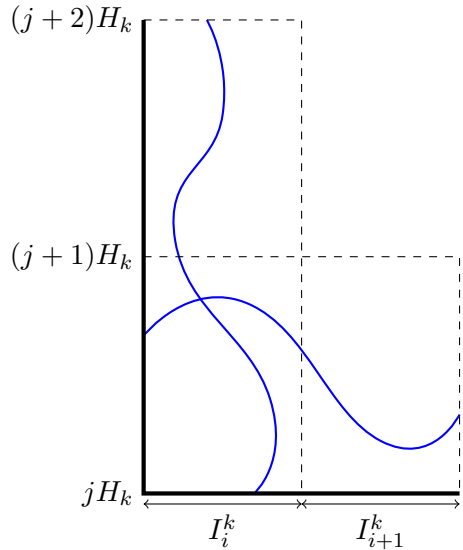
\begin{figure}[!htbp]
    \centering
    \begin{tikzpicture}[scale=.7]
    \draw[thick, blue, smooth, tension=0.9] plot coordinates {
    (3-.9,0)
    (3-.6,1.8)
    (3-2.4,4.8) 
    (3-1.5,7.2)
    (3-1.8,9)};
    \draw[thick, blue, smooth, tension=0.9] plot coordinates {
    (0,3)
    (2,3.6)
    (4.5,1)
    (6,1.5)};
    \draw[ultra thick] (0,9)--(0,0) -- (6,0);
    \draw[dashed] (0,9)--(3,9) -- (3,0);
    \draw[dashed] (0,4.5)--(6,4.5) -- (6,0);
    \draw[<->] (0,-.2)--(3,-.2);
    \draw[<->] (6,-.2)--(3,-.2);
    \node[below] at (1.5,-.2) {$I^k_i$};
    \node[below] at (4.5,-.2) {$I^k_{i+1}$};
    \node[left] at (0,0) {$jH_k$};
    \node[left] at (0,4.5) {$(j+1)H_k$};
    \node[left] at (0,9) {$(j+2)H_k$};
    \end{tikzpicture}
    \caption{An illustration of occurrence of the events $H_{i,j}^k$ and $V_{i,j}^k$ denoted by the blue open paths.}
    \label{fig_cruzamentos_hev}
\end{figure}

Also, for every $i,j,k\in\Z_+$ and $p\in(0,1)$ define the probabilities
\begin{equation*}
h_k(p; i,j)=\max_{\substack{\Lambda \; : \; I_i^k \text{ and } I_{i+1}^k\\ \text{ are good}}} \PP_p^{\Lambda}\big((H_{i,j}^k)^c\big)
\end{equation*}
and
\begin{equation*}
v_k(p;i,j)=\max_{\substack{\Lambda \; : \; I_i^k \\ \text{ is good}}} \PP_p^{\Lambda}\big((V_{i,j}^k)^c\big),
\end{equation*}
where $(H_{i,j}^k)^c$ and $(V_{i,j}^k)^c$ denote the complementary event of $H_{i,j}^k$ and $V_{i,j}^k$, respectively. Also define  
\begin{equation}
q_k(p;i.j)=\max\{h_k(p;i,j), \; v_k(p;i,j)\}.
\label{def_qk}
\end{equation}
Translation invariance allows us to write, for each $k\in\Z_+$,
\begin{equation}
q_k(p):=q_k(p;0,0)=q_k(p;i,j), \text{ for any } i,j\in\Z_+.
\label{qk}
\end{equation}

Now, we will show that, for $p$ sufficiently close to 1, the sequence $q_k(p)$ decrease fast enough in $k$. For this, we will require two auxiliary lemmas to deal with horizontal and vertical crossings, respectively.

\begin{lema}[{\bf Horizontal crossings}]
\label{lema_cota_cruz_hor}
Let $A$ given by \eqref{def_A} and let $p>\frac{1}{2}$. There exists $k_1=k_1(\alpha, \beta, \mu) \geq k_0 \in \Z_+$ such that, if
\begin{equation}
q_k(p) \leq e^{-L_k^{\left(1-\frac{\mu}{k+1}\right)}}
\label{hipotese_qk}
\end{equation}
then 
\begin{equation}
\PP_p^{\Lambda}\big((H_{0,0}^{k+1})^c\big) \leq e^{-L_{k+1}^{\left(1-\frac{\mu}{k+2}\right)}},
\label{inducao_cruz_horiz}
\end{equation}
for $k \geq k_1$ and every environment $\Lambda \in \{I_0^{k+1} \text{ is good}\} \cap \{I_1^{k+1} \text{ is good}\}$.
\end{lema}
The hypothesis $p > \frac{1}{2}$ is not important and was only adopted for convenience, to simplify the calculations.

\begin{dem}
Take an integer $k_1=k_1(\alpha, \beta, \mu) \geq k_0$ such that for $k \geq k_1$, the three inequalities below hold: 
\begin{align}
& \exp{\left[(k+2)\log A - \left(\dfrac{A^{\frac{(k+1)(k+2)}{2}}}{2^{k+1}}\right)^{1-\mu}\right]} \leq \frac{1}{32}, \label{eq_hor_k1} \\
& A^{\frac{-(2-\beta)k+3\beta-4}{2}} \leq \dfrac{1}{36}, \label{eq_hor_k1_2} \\
& \exp{\left[-\exp\left(\dfrac{A^{\frac{(k+2-\beta)(k+3)}{2}}}{2^{k+3-\beta}}\right) + A^{\frac{(k+2-\mu)(k+3)}{2}}\right]} \leq 1. \label{eq_hor_k1_3}
\end{align}

Fix an environment $\Lambda$ for which $I_0^{k+1}$ and $I_1^{k+1}$ are good intervals. By definition, both these intervals can contain at most two bad intervals at scale $k$ and, in this case, they must be adjacent.  Although the probability of crossing a bad interval on scale $k$ is small, the exponential height of the rectangles guarantees that will be a lot of attempts to do this. Indeed, let us divide the rectangle $R([0,2L_{k+1})\times[0,H_{k+1}))$ into strips of height $2H_k$ and verify whether or not we crossed these strips. For each $0\leq j\leq \left\lceil e^{L_{k+1}^{\left(1-\frac{\beta}{k+2}\right)}}\right\rceil -1$, define the events
\begin{equation*}
S_j^{k+1}=\mathcal{C}_h\big(R([0,2L_{k+1})\times[2jH_k, (2j+2)H_k))\big).
\end{equation*}
Notice if the event $H_{0,0}^{k+1}$ does not occur then none of the events $S_j^{k+1}$ can occur, that is, $\big\{ (H_{0,0}^{k+1})^c \subseteq \bigcap_j (S_j^{k+1})^c \big\}$, where $0\leq j\leq \left\lceil e^{L_{k+1}^{\left(1-\frac{\beta}{k+2}\right)}} \right\rceil -1$,  see Figure \ref{fig_evento_S}.

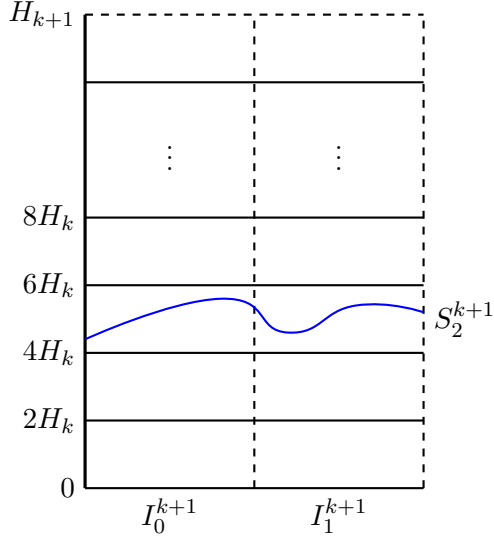
\begin{figure}[!htbp]
    \centering
    \begin{tikzpicture}[scale=.9]
      \draw[very thick] (0,0) -- (0,7);
      \draw[dashed, thick] (5,0) -- (5,7);
      \draw[dashed, thick] (0,7) -- (5,7);
      \draw[dashed, thick] (2.5,0) -- (2.5,7);
      \foreach \y in {0,1,...,4,6}
      \draw[thick] (0,\y) -- (5,\y);

      \node at (1.25,5) {$\vdots$};
      \node at (3.75,5) {$\vdots$};

      \node[left] at (0,0) {$0$};
      \node[left] at (0,1) {$2H_k$};
      \node[left] at (0,2) {$4H_k$};
      \node[left] at (0,3) {$6H_k$};
      \node[left] at (0,4) {$8H_k$};
      \node[left] at (0,7) {$H_{k+1}$};
      \node[below] at (1.25,0) {$I^{k+1}_0$};
      \node[below] at (3.75,0) {$I^{k+1}_1$};
      \node[right] at (5,2.5) {$S^{k+1}_2$};

      \draw[thick, blue, smooth, tension=0.9] plot coordinates {
         (0,2.2)
         (2,2.8)
         (3,2.3)
     (4,2.7)
     (5,2.6)};
    \end{tikzpicture}
    \caption{An illustration of the occurrence of the event $S_2^{k+1}$. Its occurrence implies the occurrence of the event $H_{0,0}^{k+1}$. }
    \label{fig_evento_S}
\end{figure}

Therefore,
\begin{equation}
\PP_p^{\Lambda}\big((H_{0,0}^{k+1})^c\big) \leq \PP_p^{\Lambda} \left(\bigcap_j (S_j^{k+1})^c \right) = \bigcap_j \PP_p^{\Lambda} \big( (S_j^{k+1})^c  \big),
\label{eq_cruz_1}
\end{equation}
with $j$ satisfying $0\leq j\leq \left\lceil e^{L_{k+1}^{\left(1-\frac{\beta}{k+2}\right)}}\right\rceil -1$. The inequality comes from inclusion of events and the equality holds since the events $(S_j^{k+1})^c$ are independent. Now, using the translation invariance of the events $S_j^{k+1}$'s in \eqref{eq_cruz_1}, we can write
\begin{equation}
\PP_p^{\Lambda}\big((H_{0,0}^{k+1})^c\big) \leq \left[\PP_p^{\Lambda}\big(S_0^{k+1}\big)^c\right]^{\left\lceil e^{L_{k+1}^{\left(1-\frac{\beta}{k+2}\right)}}\right\rceil} \leq \left[1-\PP_p^{\Lambda}\big(S_0^{k+1}\big)\right]^{\exp{\left(L_{k+1}^{\left(1-\frac{\beta}{k+2}\right)}\right)}}.
\label{probcruzhorizcomp}
\end{equation}

In order to get an upper bound for the probability of $S_0^{k+1}$, we will build a horizontal crossing within the strip $R([0,2L_{k+1})\times[0,2H_k))$ using the crossings events $H_{i,0}^k$ and $V_{i,0}^k$ in rectangles whose bases are good $k$-intervals of $I_0^{k+1}$ and $I_1^{k+1}$ while, in the rectangles whose bases are bad $k$-intervals, we will open paths at its top as follows.

Notice that the base of the strip $R([0,2L_{k+1})\times[0,2H_k))$ is divided into two parts of length $L_{k+1}$ each one. Furthermore, each of these parts have $\lfloor A^{k+2} \rfloor$ intervals of length $L_k$. Since $I_0^{k+1}$ and $I_1^{k+1}$ are good intervals, among the $2 \lfloor A^{k+2} \rfloor$ intervals in the base of the strip, we will have at most two bad intervals $L_k$, which must be adjacent, in each part of length $L_{k+1}$.

For each $l\in\{0,1\}$, denote by $j_l$ the index of the first bad $k$-interval within of $I_l^{k+1}$ and consider the interval
\begin{equation*}
I_l^* = \big( I_{j_l-1}^k \cup I_{j_l}^k \cup I_{j_l+1}^k \cup I_{j_l+2}^k \big) \cap \big( I_0^{k+1} \cup I_1^{k+1} \big) \subseteq \Z_+,
\end{equation*}
that is, $I_l^*$ is the interval formed by $I_{j_l}^k$, the $k$-interval before it and the two $k$-intervals after it (as long as are contained in $I_0^{k+1} \cup I_1^{k+1}$). Also define, for each $l\in\{0,1\}$, the path $\gamma_l^*$ formed by the edges of the form $\langle (m,2H_k), (m+1,2H_k) \rangle$, with $m \in I_l^*$, see Figure \ref{fig_cruz_faixa_horiz}. Thus, we have
\begin{equation}
\left( \bigcap_{\substack{i: I_i^k, \; I_{i+1}^k \\ \text{are goods}}} H_{i,0}^k \right) \cap \left( \bigcap_{\substack{j: I_j^k \\ \text{is good}}} V_{j,0}^k \right) \cap \{\gamma_0^* \text{ and } \gamma_1^* \text{ are open paths}\} \subseteq S_0^{k+1},
\label{evento_faixa_horiz}
\end{equation}
where $0 \leq i,j \leq 2 \lfloor A^{k+2} \rfloor -1$. 
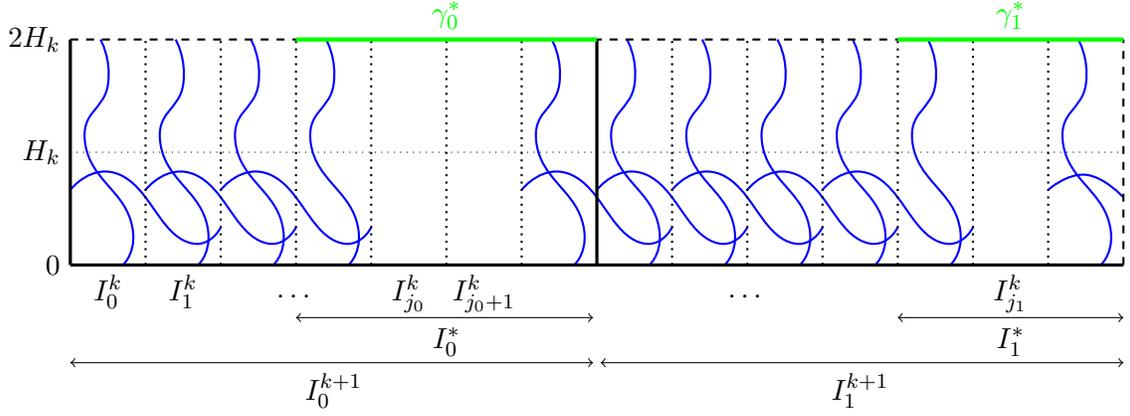
\begin{figure}[!htbp] 
    \centering
        \begin{tikzpicture}[scale=1]
     \foreach \x in {0,1,2,3,6,7,8,9,10,11,13}{
      \draw[thick, blue, smooth, tension=0.9] plot coordinates {
         (\x+.7,0)
         (\x+.8,.6)
         (\x+.2,1.6) 
         (\x+.5,2.4)
         (\x+.4,3)};
     }
 \foreach \x in {0,1,2,6,7,8,9,10}{
     \draw[thick, blue, smooth, tension=0.9] plot coordinates {
         (\x+0,1)
         (\x+.66,1.2)
         (\x+1.5,.33)
         (\x+2,.5)};
     }
\draw[thick, blue, smooth, tension=0.9] plot coordinates {
    (13+0,1)
    (13+.5,1.2)
    (13+1,.9)};
     
      \draw[very thick] (0,0) -- (0,3);
       \draw[very thick] (7,0) -- (7,3);
      \draw[very thick] (0,0) -- (14,0);
      \draw[dashed, thick] (14,0) -- (14,3);
      \draw[dashed, thick] (0,3) -- (14,3);
      \foreach \x in {1,...,13}
      \draw[thick, dotted] (\x,0) -- (\x,3);
\draw[dotted] (0,1.5) -- (14,1.5);

\draw[ultra thick,green] (3,3)--(7,3); 
\draw[ultra thick,green] (11,3)--(14,3); 
      \node[above,green] at (5,3) {$\gamma_0^*$};
      \node[above,green] at (12.5,3) {$\gamma_1^*$};
      \node[left] at (0,0) {$0$};
      \node[left] at (0,1.5) {$H_k$};
      \node[left] at (0,3) {$2H_k$};
      \node[below] at (.5,0) {$I^{k}_0$};
      \node[below] at (1.5,0) {$I^{k}_1$};
      \node[below] at (3,-0.2) {$\cdots$};
      \node[below] at (4.5,0) {$I^k_{j_0}$};
      \node[below] at (5.5,0) {$I^k_{j_0+1}$};
      \node[below] at (9,-0.2) {$\cdots$};
      \node[below] at (12.5,0) {$I^k_{j_1}$};
      
      \draw[<->](3,-.7)--(6.95,-.7);
       \node[below] at (5,-.7) {$I^*_0$};
       \draw[<->](11,-.7)--(14,-.7);
       \node[below] at (12.5,-.7) {$I^*_1$};
        \draw[<->](0,-1.3)--(6.95,-1.3);
       \node[below] at (3.5,-1.3) {$I^{k+1}_0$};
       \draw[<->](7.05,-1.3)--(14,-1.3);
       \node[below] at (10.5,-1.3) {$I^{k+1}_1$};
    \end{tikzpicture}
    \caption{An illustration of the occurrence of events $H_{i,0}^k$ and $V_{i,0}^k$ (blue open paths) and of two open paths $\gamma_0^*$ and $\gamma_1^*$ (green paths), which imply the occurrence of $S_0^{k+1}$. The intervals $I_{j_0}^k, \; I_{j_0+1}^k$ and $I_{j_1}^k$ correspond to bad $k$-intervals.}
    \label{fig_cruz_faixa_horiz}
\end{figure}
Since $I_0^*$ and $I_1^*$ are disjoint intervals, 
\begin{align}
\nonumber \PP_p^{\Lambda}(\gamma_0^* \text{ and } \gamma_1^* \text{ are open paths}) & = \PP_p^{\Lambda}(\gamma_0^* \text{ is an open path}) \PP_p^{\Lambda}(\gamma_1^* \text{ is an open path}) \\
 & \geq \big(p^{4L_k} \big)^2 = p^{8L_k}.
\label{prob_caminhos}
\end{align}
On other hand, by FKG inequality, independence of events $H$'s and $V$'s, Bernoulli's inequality, \eqref{qqc} and \eqref{hipotese_qk}, we have

\begin{align}
\nonumber \PP_p^{\Lambda} \left[\left( \bigcap_{\substack{i: I_i^k, \; I_{i+1}^k \\ \text{are goods}}} H_{i,0}^k \right)  \cap \left( \bigcap_{\substack{j: I_j^k \\ \text{is good}}} V_{j,0}^k \right) \right] &\geq \prod_{i=0}^{2 \lfloor A^{k+2} \rfloor -1} \PP_p^{\Lambda} (H_{i,0}^k) \; \PP_p^{\Lambda}(V_{i,0}^k) \\
\nonumber & \geq \big[\big(1-q_k(p)\big)^2\big]^{2A^{k+2}} \\ 
\nonumber & \geq \big(1-4 A^{k+2} \; q_k(p)\big)\\
\nonumber & \geq \big(1- 16  L_k^{\frac{2}{k+1}}\; q_k(p)\big)\\
 & \geq \left(1-16 L_k^{\frac{2}{k+1}} \; e^{- L_k^{\left(1-\frac{\mu}{k+1}\right)}} \right).
\label{prob_cruzamentos}
\end{align}
Notice that, from \eqref{cotas_Lk} and \eqref{eq_hor_k1}
\begin{align*}
L_k^{\frac{2}{k+1}} \; e^{- L_k^{\left(1-\frac{\mu}{k+1}\right)}} & = \exp{\left(\dfrac{2}{k+1}\log L_k - L_k^{1-\frac{\mu}{k+1}}\right)} \\
 & \leq \exp{\left[(k+2) \log A - \left(\dfrac{A^{\frac{(k+1)(k+2)}{2}}}{2^{k+1}} \right)^{1-\mu}\right]} \\
 & < \dfrac{1}{32}.
\end{align*}
So, from \eqref{prob_cruzamentos}, we get
\begin{equation}
\PP_p^{\Lambda} \left[\left( \bigcap_{\substack{i: I_i^k, \; I_{i+1}^k \\ \text{are goods}}} H_{i,0}^k \right)  \cap \left( \bigcap_{\substack{j: I_j^k \\ \text{is good}}} V_{j,0}^k \right) \right] \geq \dfrac{1}{2}.
\label{prob_cruz_x}
\end{equation}
From FKG inequality, \eqref{evento_faixa_horiz}, \eqref{prob_caminhos} and \eqref{prob_cruz_x} we have
\begin{equation}
\PP_p^{\Lambda}(S_0^{k+1}) \geq \frac{1}{2}\, p^{8L_k} \geq e^{-9L_k}.
\label{prob_cruz_faixa_horiz}
\end{equation}
Substituting \eqref{prob_cruz_faixa_horiz} in \eqref{probcruzhorizcomp} and using the inequality $1-x \leq e^{-x}$, we obtain
\begin{align}
\nonumber \PP_p^{\Lambda}\big((H_{0,0}^{k+1})^c\big) & \leq \left( 1- e^{-9 L_k} \right)^{\exp\left(L_{k+1}^{\left(1-\frac{\beta}{k+2}\right)}\right)} \\
& \leq \exp{\left[-\exp\left(-9 L_k + L_{k+1}^{\left(1-\frac{\beta}{k+2}\right)}\right)\right]}. 
\label{eq2}
\end{align}

Since \eqref{cotas_Lk+1_e_Lk} implies that $L_k \leq 2 L_{k+1}^{\frac{k+1}{k+3}}$, by \eqref{eq2} we get
\begin{align}
\nonumber \dfrac{\PP_p^{\Lambda}\big((H_{0,0}^{k+1})^c\big)}{\exp\left(- L_{k+1}^{\left(1-\frac{\mu}{k+2}\right)}\right)} & \leq \exp{\left[-\exp\left(-9 L_k + L_{k+1}^{\left(1-\frac{\beta}{k+2}\right)}\right) + L_{k+1}^{\left(1-\frac{\mu}{k+2}\right)} \right]} \\
\nonumber & \leq \exp{\left[-\exp\left(-18 L_{k+1}^{\left(1-\frac{2}{k+3}\right)} + L_{k+1}^{\left(1-\frac{\beta}{k+2}\right)}\right) + L_{k+1}^{\left(1-\frac{\mu}{k+2}\right)} \right]} \\
 & \leq \exp{\left[-\exp\left[ L_{k+1}^{\left(1-\frac{\beta}{k+2}\right)}\left(1-18 L_{k+1}^{\left(-\frac{2}{k+3} + \frac{\beta}{k+2}\right)} \right) \right] + L_{k+1}^{\left(1-\frac{\mu}{k+2}\right)} \right]}.
\label{eq_1}
\end{align}
Notice that, from \eqref{cotas_Lk} and \eqref{eq_hor_k1_2}, we obtain
\begin{equation}
1-18 L_{k+1}^{\left(-\frac{2}{k+3} + \frac{\beta}{k+2}\right)} = 1-18 L_{k+1}^{\left(\frac{(\beta-2)k+3\beta-4}{(k+2)(k+3)}\right)} \geq 1-18 A^{\frac{(\beta-2)k+3\beta-4}{2}} \geq \dfrac{1}{2}.
\label{eq_in}
\end{equation}
Hence, from \eqref{eq_1}, \eqref{eq_in} and \eqref{cotas_Lk},
\begin{align*}
\dfrac{\PP_p^{\Lambda}\big((H_{0,0}^{k+1})^c\big)}{\exp\left(- L_{k+1}^{\left(1-\frac{\mu}{k+2}\right)}\right)} & \leq \exp{\left[-\exp\left(\dfrac{1}{2} L_{k+1}^{\left(1-\frac{\beta}{k+2}\right)} \right) + L_{k+1}^{\left(1-\frac{\mu}{k+2}\right)} \right]} \\
& \leq \exp \left[ -\exp{\left(\dfrac{A^{\frac{(k+2-\beta)(k+3)}{2}}}{2^{k+3-\beta}}\right)} + A^{\frac{(k+2-\mu)(k+3)}{2}} \right] 
\end{align*}
and, by\eqref{eq_hor_k1_3}, the last expression is smaller than $1$ for $k\geq k_3$.
\end{dem}

So far all the estimates worked out  for any $c>0$. In order to control the estimates for the vertical crossings, we will need to restrict ourselves  to $c> 8\sqrt{\log 96}$.

\begin{lema}
\label{lemmdc}
Let $c>8\sqrt{\log 96}$ and $A(\alpha,c)$ be defined as in \eqref{def_A}. There exists $e^{-c^2/32}<\alpha_o<1$ and  $0<\mu_o<1$, such that if $(\alpha,\mu)\in (\alpha_o,1)\times (\mu_o,1)$, then
\begin{equation}
\label{eqwx1}
\big(A(\alpha,c)\big)^{\mu/2}>96.
\end{equation}.
\end{lema}

\begin{dem}
For $\alpha,\mu>0$, define
\begin{equation*}
w(\alpha,\mu)=\sqrt{32\, \log\left(\frac{96^{\frac{2}{\mu}}}{\alpha}\right)}.
\end{equation*}
Notice that $w(1,1)=\sqrt{64\log 96}=8\sqrt{\log 96}<c$ and so, by continuity of $w$ at $(1,1)$, there exists $e^{-c^2/32}<\alpha_o<1$ and $0<\mu_o<1$ such that if $(\alpha,\mu)\in(\alpha_o,1)\times(\mu_o,1)$, then $w(\alpha,\mu)<c$.  And so,
\begin{equation*}
A(\alpha,\mu)=\alpha e^{\frac{c^2}{32}}>\alpha e^{\frac{w^2(\alpha,\mu)}{32}}=96^{\frac{2}{\mu}},
\end{equation*}
which proves the lemma.
\end{dem}

\begin{lema}[{\bf Vertical Crossings}]
\label{lema_cota_cruz_vert}
Let $\alpha$, $\mu$ and $c_o$ be given as in Lemma \ref{lemmdc} and $\mu<\beta<1$. Suppose that $c\geq c_o$.  There exists a positive integer $k_2=k_2( \alpha, \beta,\mu) \geq k_0$ such that, if
\begin{equation}
q_k(p) \leq e^{-L_k^{\left(1-\frac{\mu}{k+1}\right)}}
\label{cota_qk_2}
\end{equation} 
then 
\begin{equation}
    \PP_p^{\Lambda}\big((V_{0,0}^{k+1})^c\big) \leq e^{-L_{k+1}^{\left(1-\frac{\mu}{k+2}\right)}}
\label{inducao_cruz_vert}
\end{equation}
for $k \geq k_2$ and every environment $\Lambda \in \{I_0^{k+1} \text{ is good}\}$.
\end{lema}

\begin{dem}
Take $k$ sufficiently large such that the two following conditions hold:
\begin{equation}
L_k^{\left(1-\frac{\mu}{k+1}\right)} > 19 \ln 2,
\label{cond_k1w}
\end{equation}
\begin{equation}
e^{\frac{-L_k^{\left(1-\frac{\mu}{k+1}\right)}+19 \ln 2 }{6}} < \dfrac{1}{2}.
\label{cond_k2w}
\end{equation}

Fix an environment $\Lambda$ for which $I_0^{k+1}$ is a good interval. By definition, this interval can contain at most two bad intervals on the scale $k$ and, in this case they must be adjacent. In this way, either each $I_i^k$ is good for every $i=0,1,\dots, \Big\lfloor \frac{A^{k+2}}{2}\Big\rfloor -2$ or $I_i^k$ is good for every $i= \Big\lfloor \frac{A^{k+2}}{2} \Big\rfloor +1, \Big\lfloor \frac{A^{k+2}}{2} \Big\rfloor + 2, \dots, L_k^{\frac{2}{k+1}} - 1$. Assume, without loss of generality, that the first case holds, and let
\begin{equation}
M_k = \Bigg\lfloor\dfrac{A^{k+2}}{2} \Bigg\rfloor -1,
\label{rede_reescalada}
\end{equation}
namely, $M_k$ represents the number of good intervals $I_i^k$ in $I_0^{k+1}$.

In order to estimate the probability of the event $V_{0,0}^{k+1}$ we will consider a following rescaled lattice: each rectangle $R\big(I_i^k \times [jH_k,(j+1)H_k)\big)$ will correspond to a vertex $(i,j)$, for all $i,j \in \Z_+$. Such vertex $(i,j)$ is declared open if the event $H_{i,j}^k\cap V_{i,j}^k$ occurs in the original lattice. Notice that the resulting percolation process in this renormalized lattice is a dependent process, since the state of one vertex $(i,j)$ depend on the state of other six vertices so that there is 6 vertices $(i',j')$ such that the events $\{(i,j) \text{ is open}\}$ and $\{(i',j') \text{ is open}\}$ are dependents, see Figure \ref{fig_vertices_dependencia} for more details. 
\begin{figure}[H]
    \centering
 \begin{tikzpicture}[scale=.7]        
        \foreach \x in {0,...,5}
        \draw[dotted] (\x,0) -- (\x,8);
        \foreach \y in {0,1.5,...,7.5}
        \draw[dotted] (0,\y) -- (5.5,\y);
        
         \draw[thick, blue, smooth, tension=0.9] plot coordinates {
           (2+.7,3)
           (2+.8,3.6)
           (2+.2,4.6) 
           (2+.5,5.4)
           (2+.4,6)};
  
\draw[thick, blue, smooth, tension=0.9] plot coordinates {
    (2+0,4)
    (2+.66,4.2)
    (2+1.5,3.33)
    (2+2,3.5)};

        \draw[ultra thick] (0,0) -- (5.5,0);
        \draw[ultra thick] (0,0) -- (0,8);
        \draw[ultra thick] (2,3) rectangle (3,6);
         \draw[ultra thick] (2,3) rectangle (4,4.5);
        \node[below] at (2.5,0) {$I^k_i$};
          \node[below] at (1,-.2) {$\cdots$};
            \node[below] at (5,-.2) {$\cdots$};
         \node[below] at (3.5,0) {$I^k_{i+1}$};
        
        \node[left] at (0,3) {$jH_k$};
        \node[left] at (0,4.5) {$(j+1)H_k$};
        \node[left] at (0,6) {$(j+2)H_k$};
        
        \draw[->, ultra thick] (6,3.5) -- (7,3.5);
        
      \foreach \x in {0,...,5}
      \draw[dotted] (8+\x,0) -- (8+\x,8);
      \foreach \y in {0,1.5,...,7.5}
      \draw[dotted] (0+8,\y) -- (5.5+8,\y);
      
         \draw[ultra thick] (8,0) -- (8+5.5,0);
      \draw[ultra thick] (8,0) -- (8,8);
            \foreach \x/\y in {2.5/2.25, 3.5/2.25, 3.5/3.75, 1.5/3.75, 1.5/5.25, 2.5/5.25}
            \filldraw[blue] (8+\x,\y) circle (6pt);
            
            \filldraw[black] (10.5,3.75) circle (6pt);
            
            \node[below] at (10.5,0) {$i$};
            \node[left] at (8,3.75) {$j$};
    \end{tikzpicture}
    \caption{The figures on the left and on the right illustrate, respectively, the occurrence of event $H_{i,j}^k\cap V_{i,j}^k$ and the renormalized lattice. The vertex $(i,j)$ is represented by the black ball and its state depends on the states of the six vertices represented by the blue balls.}
    \label{fig_vertices_dependencia}
\end{figure}
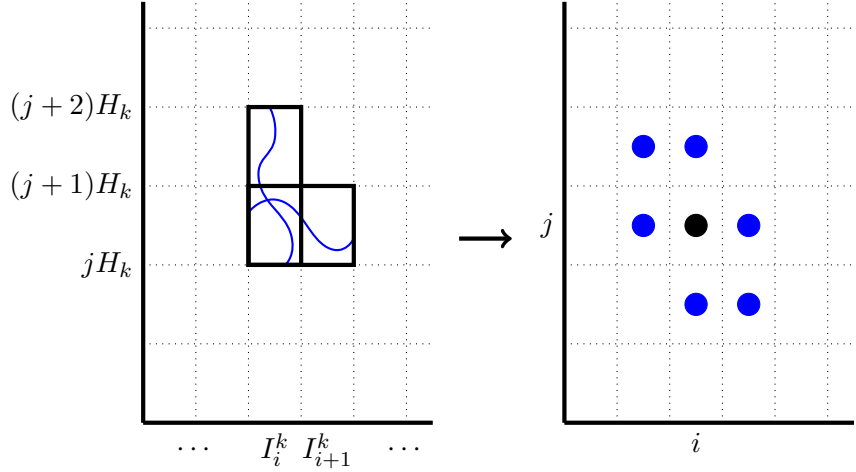

Consider the rectangle
\begin{equation*}
R_0^{k+1}=R\left(\left[0, M_k \right) \times \left[0, 4 \bigg\lceil e^{L_{k+1}^{\left(1-\frac{\beta}{k+2}\right)}} \bigg\rceil \right)\right),
\end{equation*}
where the right side of the above equality is given by \eqref{retangulo} and $M_k$ is given by \eqref{rede_reescalada}, and let its vertical crossing event be denoted by $\mathcal{C}_v(R_0^{k+1})$. Notice that
\begin{equation}
\PP_p^{\Lambda}(V_{0,0}^{k+1}) \geq \PP\big(\mathcal{C}_v(R_0^{k+1})\big).
\label{cota_prob_cruz_vert}
\end{equation}

Suppose that the the event $\mathcal{C}_v(R_0^{k+1})$ does not occur. Then there must be a sequence of distinct vertices, namely $(i_0,j_0), (i_1,j_1), \dots , (i_n,j_n)$ in $R_0^{k+1}$ such that the following three conditions hold:
\begin{enumerate}[(i)]
    \item $\displaystyle\max_{1 \leq l \leq n}\big\{ |i_l - i_{l-1}|, |j_l - j_{l-1}|\big\} = 1$,
    \item $(i_0,j_0)\in\{0\} \times \left[0, 4 \bigg\lceil e^{L_{k+1}^{\left(1-\frac{\beta}{k+2}\right)}} \bigg\rceil \right]$ and $(i_n,j_n)\in\{M_k\} \times \left[0, 4 \bigg\lceil e^{L_{k+1}^{\left(1-\frac{\beta}{k+2}\right)}} \bigg\rceil \right]$, and
    \item $(i_l,j_l)$ is closed for every $l=0,1,\dots,n$.
\end{enumerate}

Notice that there are at most $4 \bigg\lceil e^{L_{k+1}^{\left(1-\frac{\beta}{k+2}\right)}}\bigg\rceil 8^n$ such sequences with $n+1$ vertices that satisfies the conditions $(i)$ and $(ii)$ above. Moreover, the probability of a vertex in $R_0^{k+1}$ be declared closed is at most $2q_k(p)$. Also, by the dependence in the rescaled lattice, for every set with $n+1$ vertices, there exists at least $\frac{n}{6}$ vertices that have been declared open or not, independently of each other. Then,
\begin{align}
\nonumber \PP_p^{\Lambda}\big(\mathcal{C}_v(R_0^{k+1})^c\big) & \leq \PP\left(\begin{array}{c}\text{there is a sequence of } n+1 \text{ vertices in } R_0^{k+1} \\  \text{ satisfying the conditions (i), (ii) and (iii)}\end{array}\right) \\
\nonumber & \leq \sum_{n \geq M_k}  4 \bigg\lceil e^{L_{k+1}^{\left(1-\frac{\beta}{k+2}\right)}} \bigg\rceil 8^n \big(2 q_k(p) \big)^{\frac{n}{6}} \\
\nonumber & \leq  4 \bigg\lceil e^{L_{k+1}^{\left(1-\frac{\beta}{k+2}\right)}} \bigg\rceil \sum_{n \geq M_k}  8^n \left(2 \; e^{-L_k^{\left(1-\frac{\mu}{k+1}\right)}} \right)^{\frac{n}{6}} \\
\nonumber & =  4 \bigg\lceil e^{L_{k+1}^{\left(1-\frac{\beta}{k+2}\right)}} \bigg\rceil \sum_{n \geq M_k}  2^{\frac{19 n}{6}} \exp{\left(-\dfrac{n}{6} \; L_k^{\left(1-\frac{\mu}{k+1}\right)} \right)} \\
 & \leq 8 \bigg\lceil e^{L_{k+1}^{\left(1-\frac{\beta}{k+2}\right)}} \bigg\rceil 2^{\frac{19 M_k}{6}} \exp{\left(-\dfrac{M_k}{6}\; L_k^{\left(1-\frac{\mu}{k+1}\right)} \right)},
 \label{eq_3}
\end{align}
where the third inequality follows from \eqref{cota_qk_2} and the last inequality follows from \eqref{cond_k1w} and \eqref{cond_k2w}. Thus,
\begin{align}
\nonumber \PP_p^{\Lambda}\big((V_{0,0}^{k+1})^c\big) & \leq 8 \; e^{L_{k+1}^{\left(1-\frac{\beta}{k+2}\right)}} 2^{\frac{19}{6}\left(\frac{\lfloor A^{k+2} \rfloor}{2}-1\right)} \exp{\left[-\dfrac{1}{6} \left(\dfrac{\lfloor A^{k+2}\rfloor}{2}-1\right) \; L_k^{\left(1-\frac{\mu}{k+1}\right)} \right]} \\
\nonumber & \leq 8 \cdot 2^{\frac{19}{3} L_k^{\frac{2}{k+2}}} \exp{\left(L_{k+1}^{\left(1-\frac{\beta}{k+2}\right)} + \dfrac{L_k^{\left(1-\frac{\mu}{k+1}\right)}}{6} - \dfrac{L_k^{\frac{2}{k+1}} L_k^{\left(1-\frac{\mu}{k+1}\right)}}{24}\right)} \\
 & \leq 8 \exp{\left(\frac{19 \log 2}{3} L_{k+1}^{\frac{2}{k+2}}+L_{k+1}^{\left(1-\frac{\beta}{k+2}\right)} + \dfrac{2L_{k+1}^{\left(1-\frac{(2+\mu)}{k+3}\right)}}{6} - \dfrac{\frac{1}{4}L_{k+1}^{\left(1-\frac{\mu}{k+3}\right)}}{24}\right)}, 
 \label{eq_4}
\end{align} 
where the first inequality follows from \eqref{cota_prob_cruz_vert}, \eqref{eq_3} and \eqref{rede_reescalada}, the second inequality follows from \eqref{cotas_chao_Ak} and \eqref{qqc}, and the third inequality follows from \eqref{cotas_Lk+1_e_Lk}. So, we have 
\begin{align}
\nonumber & \dfrac{\PP_p^{\Lambda}\big((V_{0,0}^{k+1})^c\big)}{\exp\left(-L_{k+1}^{\left(1-\frac{\mu}{k+2}\right)}\right)} \\
\nonumber & \leq 8 \exp{\left(\frac{19 \log 2}{3} L_{k+1}^{\frac{2}{k+2}}+L_{k+1}^{\left(1-\frac{\beta}{k+2}\right)} + \dfrac{L_{k+1}^{\left(1-\frac{(2+\mu)}{k+3}\right)}}{3} - \dfrac{L_{k+1}^{\left(1-\frac{\mu}{k+3}\right)}}{96} +L_{k+1}^{\left(1-\frac{\mu}{k+2}\right)}\right)}  \\
\nonumber & = 8 \exp{\left[-L_{k+1}^{\left(1-\frac{\mu}{k+3}\right)} \left(\frac{1}{96} - \frac{L_{k+1}^{-\frac{2}{k+3}}}{3} - L_{k+1}^{\frac{(\mu-\beta)k+2\mu-3\beta}{(k+2)(k+3)}} - L_{k+1}^{\frac{-\mu}{(k+2)(k+3)}}-\frac{19  \log 2}{3}L_{k+1}^{\frac{-k^2-3k+(k+2)\mu}{(k+2)(k+3)}}\right)\right]}\\
\nonumber & \leq 8 \exp{\left[-L_{k+1}^{\left(1-\frac{\mu}{k+3}\right)} \left(\frac{1}{96} - \frac{1}{3}\, \frac{A^{-(k+2)}}{2^{-\frac{2(k+2)}{k+3}}} - \frac{A^{\frac{-(\beta-\mu)k+2\mu-3\beta}{2}}}{2^{\frac{-(\beta -\mu)k+2\mu-3\beta}{k+3}}} - \frac{A^{-\frac{\mu}{2}}}{2^{-\frac{\mu}{k+3}}}-\frac{19 \,\, 2^k \,\log 2}{ A^{\frac{k^2}{2}}} \right)\right]} \\
 & := 8 \exp{\left(-L_{k+1}^{\left(1-\frac{\mu}{k+3}\right)} \, B(k, \mu, \beta, c) \right)},
\label{eq_wx1}
\end{align}
where the first and the second inequalities follow from \eqref{eq_4} and \eqref{cotas_Ax}, respectively.
By Lemma \ref{lemmdc} we have $A^{\frac{\mu}{2}} > 96$, then
\begin{equation*}
\lim_{k \rightarrow \infty} B(k, \mu, \beta, c) = \dfrac{1}{96} - \dfrac{1}{A^{\frac{\mu}{2}}} > 0.   
\end{equation*}
So,  for $k$ sufficiently large, we have $B(k, \mu, \beta, c)>\frac{1}{2}\left( \frac{1}{96}-\frac{1}{A^{\mu/2}}\right)$. Therefore, there is an integer $k_2=k_2(\alpha, \beta, \mu) \geq k_0$ such that for $k \geq k_2$ we have \eqref{cond_k1w}, \eqref{cond_k2w} and
\begin{equation*}
8 \exp{\left(-L_{k+1}^{\left(1-\frac{\mu}{k+3}\right)} \, B(k, \mu, \beta, c) \right)} \leq 1,
\end{equation*}
which concludes the proof.
\end{dem}

The proposition bellow gives us the decay of $q_k$.
\begin{prop}
\label{lema_cota_qk} 
Let $\alpha$, $\mu$ and $c_o$ be given as in Lemma \ref{lemmdc}. Suppose that $c\geq c_o$. Then, there exists $k_3=k_3( \alpha, \beta, \mu) \in \Z_+$ and $p_0=p_0(c, \alpha, \beta,\mu)<1$ such that, for all $k \geq k_3$ and $p\geq p_0$, we have
\begin{equation}
q_k(p) \leq e^{-L_k^{\left(1-\frac{\mu}{k+1}\right)}}.
\label{inducao_qk}
\end{equation}
\end{prop}

\begin{dem} The proof of will be done by induction on $k$.
Let $k_3 = \max\{k_1, k_2\} \geq k_0$, where $k_1$ and $k_2$ are given by Lemmas \ref{lema_cota_cruz_hor} and \ref{lema_cota_cruz_vert}, respectively. Since $q_{k_3}(p)$ goes to $0$ as $p$ goes to $1$, then there is $p_0=p_0(c,\alpha,\beta,\mu)<1$ such that, for all $p \geq p_0$,
\begin{equation*}
q_{k_3}(p) \leq \exp{\left[- L_{k_3}^{1-\mu}\right]}.
\end{equation*}
Lemmas \ref{lema_cota_cruz_hor} and \ref{lema_cota_cruz_vert} imply, respectively, that $\PP_p^{\Lambda}\big((H_{0,0}^{k+1})^c\big) \leq \exp{\left[-L_{k+1}^{\left(1-\frac{\mu}{k+2}\right)}\right]}$ and $\PP_p^{\Lambda}\big((V_{0,0}^{k+1})^c\big) \leq \exp{\left[-L_{k+1}^{\left(1-\frac{\mu}{k+2}\right)}\right]},$ for all $k \geq k_3$. Therefore, from \eqref{def_qk}, we conclude that
$q_k(p) \leq \exp{\left[-L_k^{\left(1-\frac{\mu}{k+1}\right)}\right]}$, for all $k \geq k_3$.
\end{dem}

The next result is the last one we will need to prove the Theorem \ref{teo_trans_fase1}.
\begin{lema} \label{lemfinz} 
Let  $0<\mu<1$ and $A>1$ given as in \eqref{def_A}. Then there is a positive integer $k_4=k_4(\alpha,\mu)$, such that for all $k\geq k_4$, we have 
\begin{equation}
\label{eqfinx}\sum_{k\geq k_4} A^{k+2}\exp{\left[-\left(\frac{A^{\frac{(k+1)(k+2)}{2}}}{2^{k+1}}\right)^{1-\frac{\mu}{k+1}}\right]} < \frac{1}{2}.
\end{equation}
\end{lema}

\begin{dem} Let $k_4 = k_4(\alpha,\mu)$ be a positive integer such that $k_4>\max\left\{4\sqrt{2},\frac{4\log 2}{\log A}\right\}$ and for all $k\geq k_4$ we have 
\begin{equation}\label{eqxx5} 
\frac{A^{k+2}e^{-\frac{A^{\frac{k^2}{8}}}{\sqrt{2}}}}{1-e^{-\frac{A^{\frac{k^2\log A}{8}}}{8\sqrt{2}}}}<\frac{1}{2}.
\end{equation}
For any $k\geq k_4$, since $k>\dfrac{4\log 2}{\log A}$, then 
\begin{equation}\label{eqyy1}
\dfrac{A^{\frac{(k+1)(k+2)}{2}}}{2^{k+1}}\geq \dfrac{A^{\frac{k^2}{2}}}{2^{k+1}} =\frac{A^{\frac{k^2}{4}}}{2} \left( \frac{A^{\frac{k}{4}}}{2}\right)^k> \dfrac{A^{\frac{k^2}{4}}}{2}.
\end{equation}  
Moreover, $\mu<1$ implies $1-\dfrac{\mu}{k+1}>\dfrac{1}{2},$ for all $k\geq 1$. Also, since $A>1$, then 
\begin{equation}
A^{\frac{l^2+2k_4l}{8}}-1=e^{\dfrac{l^2+2k_4l}{8} \log A}-1\geq \dfrac{l^2+2k_4l}{8} \log A=\dfrac{l+2k_4}{8}\, l\,  \log A, 
\label{eq_5}
\end{equation}
for all $l \geq 1$. 
Therefore, we have 
\begin{eqnarray*}\sum_{k\geq k_4} A^{k+2}e^{-\left(\frac{A^{\frac{(k+1)(k+2)}{2}}}{2^{k+1}}\right)^{1-\frac{\mu}{k+1}}}&\leq& A^2\sum_{k\geq k_4}A^k e^{-\frac{A^{\frac{k^2}{8}}}{\sqrt{2}}}\\
 \mbox{} & =  & A^{k_4+2}e^{-\frac{A^{\frac{k_4^2}{8}}}{\sqrt{2}}}\left( 1+\sum_{l\geq 1} A^l e^{-A^{\frac{k_4^2}{8}}\left(\frac{A^{\frac{l^2+2k_4l}{8}} -1}{\sqrt{2}}\right)} \right)  
 \\
 \mbox{} & \leq   & A^{k_4+2}e^{-\frac{A^{\frac{k_4^2}{8}}}{\sqrt{2}}}\left( 1+\sum_{l\geq 1}  e^{-\left(A^{\frac{k_4^2}{8}}\frac{(l+2k_4)}{8\sqrt{2}}-1\right)\, l\, \log A} \right)
 \\
 \mbox{} & \leq   & A^{k_4+2}e^{-\frac{A^{\frac{k_4^2}{8}}}{\sqrt{2}}}\left( 1+\sum_{l\geq 1}  e^{-\frac{ A^{\frac{k_4^2}{8}}\, \log A}{8\sqrt{2}}\, l^2} \right)
 \\
 \mbox{} & \leq   & A^{k_4+2}e^{-\frac{A^{\frac{k_4^2}{8}}}{\sqrt{2}}}\left( 1+\sum_{l\geq 1}  e^{-\frac{ A^{\frac{k_4^2}{8}}\, \log A}{8\sqrt{2}}\, l} \right)\\
 \mbox{} & = & \frac{A^{k_4+2}e^{-\frac{A^{\frac{k_4^2}{8}}}{\sqrt{2}}}}{1-e^{-\frac{A^{\frac{k_4^2\log A}{8}}}{8\sqrt{2}}}},
\end{eqnarray*} 
where the first, the second and the third inequalities follow from \eqref{eqyy1}, \eqref{eq_5} and the fact that $k_4> 4\sqrt{2}$, respectively. From the last equality above and \eqref{eqxx5} we conclude the proof of the lemma.
\end{dem}

\section{Proof of Theorem \ref{teo_trans_fase1}} \mbox{}\\
\label{secao_prova_teo}

Now, we are in condition to prove Theorem \ref{teo_trans_fase1}, which gives us the non-trivial phase transition for our model. The absence of percolation for small values of $p$, given by \eqref{eq_ausencia_perc} in Theorem \ref{teo_trans_fase1}, has been proved in \cite{marcos}. Using the results obtained in Sections \ref{secao_ambiente} and \ref{secao_cruzamentos}, we will show that there is a sufficiently large $p<1$ such that, for almost all realizations of $\Lambda$, the event $\{(0,0)\leftrightarrow\infty\}$ occurs with strictly positive probability.

The idea of the proof is to show percolation in the original model, indirectly, by considering a new model, which percolates. This model will be given by an aperiodic random variable $\widetilde{\xi}$ taking integer values greater than $\widetilde{c}$, given by \eqref{cota_c_teo}, such that it will be dominated by the original one. 

Let $\xi$ be any positive random variable such that $\E(\xi e^{c(\log \xi)^{1/2}} \mathbb{1}_{\{\xi \geq 1\}}) < \infty$. Let $a_1$ and $a_2$ be two positive integers such that $a_1 < a_2$ and $\PP(\lceil \xi \rceil = a_1)>0$ and $\PP(\lceil \xi \rceil = a_2)>0$ and consider $a_0 = \max \{a_2, \lceil\, \widetilde{c}\, \rceil\} + 1$. Define a new random variable $\widetilde{\xi}$ by
\begin{equation}
\widetilde{\xi} = \begin{cases}
\lceil \xi \rceil, \text{ if } \xi > a_0 \\
a_0, \text{ if } \xi \in (a_1, a_0] \\
a_0 - 1, \text{ if } \xi \leq a_1
\end{cases}.
\label{def_xi_til}
\end{equation}

Notice that $\widetilde{\xi}$ is an aperiodic, integer-valued random variable taking values greater than $\widetilde{c}$ and $\widetilde{\xi} \geq \xi$. Moreover, $\E(\widetilde{\xi} e^{c(\log \widetilde{\xi})^{1/2}}) < \infty$, since $\widetilde{\xi} e^{c(\log \widetilde{\xi})^{1/2}} \leq 2^{c+1} \xi e^{c(\log \xi)^{1/2}}$, for all $\xi \geq a_0$.

Suppose that Theorem \ref{teo_trans_fase1} is valid for $\widetilde{\xi}$. So, by stochastic domination, the same will be true for the original model.

Next we will prove Theorem \ref{teo_trans_fase1} for the model given by the variable $\widetilde{\xi}$. For the sake of notation we will replace $\widetilde{\xi}$ with $\xi$. Being $\xi$ an aperiodic and integer-valued random variable taking values greater than $\widetilde{c}$, we can apply all the results previously obtained in Section \ref{cap3}.

Let $\alpha$, $\mu$ and $c_o$ be given as in Lemma \ref{lemmdc} and $\mu<\beta<1$. Suppose that $c\geq c_o$. Notice that \eqref{prob_blocos_bons_new} means that with strictly positive probability there is an environment $\Lambda$ such that the intervals $I_0^k$ and $I_1^k$ are good, as well as their subintervals, for all scale $k \geq k_0$, where $k_0$ is given in Lemma \ref{lema_cota_pk_lindvall}. So let $\Lambda$ be fix an such environment.

\begin{figure}[!htbp]
    \centering
       \begin{tikzpicture}[scale=2]          
        \filldraw[gray] (0,0) rectangle(.2,.1);
      \foreach \x in {0,.2,...,.8} { \draw[dotted] (\x,0) -- (\x,.5);
        \draw[thin, blue, smooth, tension=0.9] plot coordinates {
          (\x+.14,0)
          (\x+.16,.1)
          (\x+.04,.3) 
          (\x+.1,.4)
          (\x+.08,.5)};
      \draw[thin, blue, smooth, tension=0.9] plot coordinates {
          (\x+0,.2)
          (\x+.133,.24)
          (\x+.3,.066)
          (\x+.4,.1)};
  }
       \foreach \x in {0.25,.5}
           \draw[dotted] (0,\x) -- (1,\x);
 \draw[dashed] (3,0) -- (3,3.3);
       \foreach \x in {0,...,2}{
       \draw[dashed] (\x,0) -- (\x,3.3);
         \draw[thick, blue, smooth, tension=0.9] plot coordinates {
           (\x+.7,0)
           (\x+.8,.6)
           (\x+.2,1.6) 
           (\x+.5,2.4)
           (\x+.4,3)};
    }
       \foreach \y in {0,1.5,3}
       \draw[dashed] (0,\y) -- (3.5,\y);
      \foreach \y in {0,1} { \draw[thick, blue, smooth, tension=0.9] plot coordinates {
           (\y+0,1)
           (\y+.66,1.2)
           (\y+1.5,.33)
           (\y+2,.5)};
       }
       \draw[thick, blue,smooth, tension=0.9] plot coordinates {
           (2+0,1)
           (2+.66,1.2)
           (2+1.3,.6)};
   
       \draw[ultra thick] (0,0) -- (3.5,0);
       \draw[ultra thick] (0,0) -- (0,3.3);
       
       \draw (-.2,0)--(-.2,3);
       \foreach \x in {0,.1,.25,.5,1.5,3}{
       \draw (-.25,\x)--(-.15,\x);}
       
        \draw (0,-.2)--(3,-.2);
       \foreach \x in {0,.2,...,1,2,3}
       \draw (\x,-.25)--(\x,-.15);
       \node[below] at (3,-.25) {$3L_k$};
      \node[below] at (2,-.25) {$2L_k$};
    \node[below] at (1,-.25) {$L_k$};
      \node[below] at (.3,-.2) {\tiny $L_{k-1}$};
        \node[below] at (3,-.25) {$3L_k$};
        \node[left] at (-.2,1.5) {$H_k$};
                   \node[left] at (-.2,3) {$2H_k$};
                   \node[left] at (-.2,.5) {\tiny $2H_{k-1}$};
                        \node[left] at (-.2,.15) {\tiny $2H_{k-3}$};
                        \node[left] at (-.2,0) { $0$};
                          \node[below] at (0,-.2) { $0$};
   \end{tikzpicture}
    \caption{An illustration of the intersections between the events $H_{i,0}^k$ and $V_{i,0}^k$ for different scales,  showing the existence of the infinite cluster.}
    \label{fig_cluster_infinito}
\end{figure}
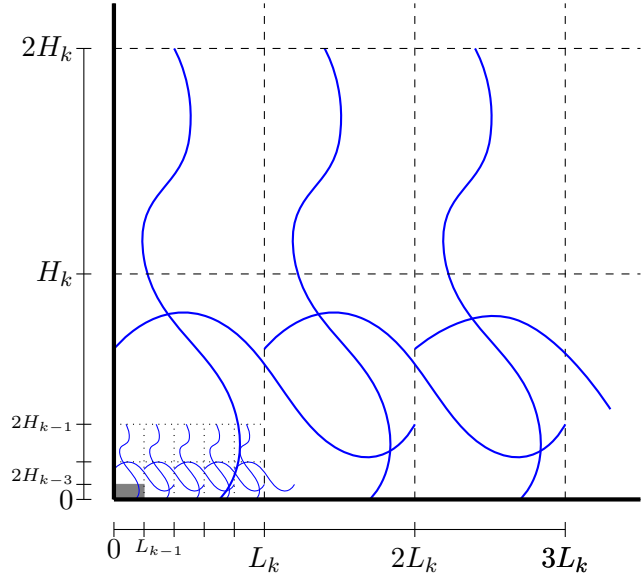

Let $k_5=k_5(\alpha, \beta, \mu)=\max\{k_3,k_4\}$, where $k_3$ and $k_4$ be given as in Proposition \ref{lema_cota_qk} and Lemma \ref{lemfinz}, respectively.  According to definition of the events $H_{i,j}^k$ and $V_{i,j}^k$ we have that 
\begin{equation}
\bigcap_{k\geq k_5} \bigcap_{i=0}^{\Big\lfloor \frac{A^{k+2}}{2}\Big\rfloor -2} (H_{i,0}^k \cap V_{i,0}^k) \subseteq \{\text{there is an infinite open cluster}\},
\label{cluster_infinito}
\end{equation}
see Figure \ref{fig_cluster_infinito}.
 
By Proposition \ref{lema_cota_qk} there is $p$ sufficiently close to $1$, independently of $\Lambda$, such that for all $k \geq k_3$, we have $q_k(p) \leq e^{-L_k^{\left(1-\frac{\mu}{k+1}\right)}}$. So, since $k_5 > k_3$, by FKG inequality, \eqref{qk}, Bernoulli's inequality, \eqref{inducao_qk} and \eqref{cotas_Lk}, we have
\begin{align}
\nonumber \PP_{p}^{\Lambda}\left(\bigcap_{k\geq k_5} \bigcap_{i=0}^{\Big\lfloor \frac{A^{k+2}}{2}\Big\rfloor -2} (H_{i,0}^k \cap V_{i,0}^k)\right)  & \geq \prod_{k\geq k_5} \prod_{i=0}^{\Big\lfloor \frac{A^{k+2}}{2}\Big\rfloor -2} \PP_{p}^{\Lambda} (H_{i,0}^k) \PP_{p}^{\Lambda}(V_{i,0}^k) \\
\nonumber & \geq \prod_{k\geq k_5} \Big[\big(1 - q_k(p) \big)^2\Big]^{\Big\lfloor \frac{A^{k+2}}{2}\Big\rfloor -2} \\
\nonumber & \geq \prod_{k \geq k_5} \big(1-q_k(p)\big)^{A^{k+2}} \\
\nonumber & \geq \prod_{k \geq k_5} \big(1 - A^{k+2} q_k(p)\big) \\
\nonumber & \geq \prod_{k \geq k_5} \left[1 - A^{k+2} \exp{\left(-L_k^{\left(1-\frac{\mu}{k+1}\right)}\right)} \right] \\
\nonumber & \geq 1 - \sum_{k\geq k_5} A^{k+2} \exp{\left(-L_k^{\left(1-\frac{\mu}{k+1}\right)}\right)} \\
\nonumber & \geq 1 - \sum_{k \geq k_5} A^{k+2} \exp{\left[-\left(\dfrac{A^{\frac{(k+1)(k+2)}{2}}}{2^{k+1}}\right)^{1-\frac{\mu}{k+1}}\right]}\\
 & > \frac{1}{2},
\label{prob_cluster}
\end{align}
where in the sixth and in the last inequalities we have used \eqref{eqfinx}, and this concludes the proof of the theorem. \hfill \qed

\section*{Acknowledgements} 
I.G. was partially supported by CAPES,  
P.C.L. was partially supported by CNPq, M. S. was partially supported by Instituto Serrapilheira (grant 5812), R. S. was partially supported by CNPq, CAPES and FAPEMIG  (grants APQ-00868-21 and RED-00133-21).

\vspace{1cm}

\begin{tabular}{cl}
 \hspace{-1cm}$^*$   &\hspace{-0.9cm}  \textsc{Universidade Federal de Minas Gerais (UFMG),}\\
  &\\
  \hspace{-1cm}$^\dagger$   &\hspace{-0.9cm} \textsc{Instituto Federal de Minas Gerais (IFMG) - Campus Conselheiro Lafaiete}\\
  
  &\\
  \hspace{-1cm}$^\diamond$   &\hspace{-0.9cm}  \textsc{Universidade Federal da Bahia (UFBA),}
 
\end{tabular}

\vspace{0.5cm}

\begin{tabular}{ll}
\hspace{-0.48cm} E-mails:    & \texttt{isadora.guedes@ifmg.edu.br},  \\
     & \texttt{lima@mat.ufmg.br}, \\
     & \texttt{marcospy6@gmail.com},\\
     &
     \texttt{rsanchis@mat.ufmg.br}.
\end{tabular}

\end{document}